\documentclass[11pt]{article}
\usepackage{amssymb, amsmath,  amsfonts,dsfont, mathrsfs,euscript,eufrak}
\usepackage{textcomp}
\usepackage{latexsym}
\usepackage{indentfirst}
\usepackage[normalem]{ulem}

\usepackage{graphicx}
\usepackage[colorlinks=true, allcolors=blue]{hyperref}

\newtheorem{Theorem}{Theorem}
\newtheorem{Lemma}{Lemma}
\newtheorem{Proposition}{Proposition}

\newtheorem{Remark}{Remark}

\newcommand{\N}{\mathbb N}
\newcommand{\R}{\mathbb R}
\newcommand{\Z}{\mathbb Z}
\newcommand{\Q}{\mathbb{Q}}
\newcommand{\CC}{\mathbb C}

\newcommand{\RID}{\boldsymbol{Q}}
\newcommand{\ID}{\boldsymbol{I}}

\newcommand{\lng}{\langle}
\newcommand{\rng}{\rangle} 

\newcommand{\id}{\mathds{1}}
\newcommand{\idd}{\mathbb{I}}
\newcommand{\dd}{d}
\newcommand{\Ln}{\mathrm{Ln}\,}
\newcommand{\Arg}{\mathrm{Arg}\,}

\newcommand{\e}{\varepsilon}

\newcommand{\poly}{\tilde{g}}
\newcommand{\sgn}{\mathrm{sgn}}
\newcommand{\ind}{\mathfrak{m}}
\newcommand{\BV}{\boldsymbol{V}}
\newcommand{\BVC}{\boldsymbol{V}_{\!\!\CC}}
\newcommand{\tphi}{\hat{\varphi}}
\newcommand{\subs}{\mathfrak{u}}
\newcommand{\X}{\mathcal{X}}

\newcommand{\tr}{\bar{r}}
\newcommand{\tbeta}{\bar{\beta}}

\newcommand{\ImagPart}{\mathrm{Im}}

\begin{document}

\title{The class $\boldsymbol{Q}$ and mixture distributions with dominated continuous singular parts}

\author{A. A. Khartov$^{1,2,3,}$\footnote{Email addresses: \texttt{khartov.a@iitp.ru}, \texttt{alexeykhartov@gmail.com}} }

\footnotetext[1]{Institute for Information Transmission Problems (Kharkevich Institute) of Russian Academy of Sciences, Bolshoy Karetny per. 19, build.1, 127051 Moscow, Russia.}
\footnotetext[2]{Smolensk State University, 4 Przhevalsky st., 214000 Smolensk, Russia. }
\footnotetext[3]{ITMO University, 49 Kronverksky Pr., 197101 Saint-Petersburg, Russia.}

\maketitle
\begin{abstract}
	We consider a new class $\boldsymbol{Q}$ of distribution functions $F$ that have the property of rational-infinite divisibility: there exist some infinitely divisible distribution functions $F_1$ and $F_2$ such that $F_1=F*F_2$. A distribution function of the class $\boldsymbol{Q}$  is quasi-infinitely divisible in the sense that its characteristic function admits  the L\'evy-type representation with a ``signed spectral measure''. This class   is a wide natural extension of the fundamental class of  infinitely divisible distribution functions and it is actively studied now. We are interested in conditions for a distribution function $F$ to belong  to the class $\boldsymbol{Q}$ for  the unexplored case, where $F$ may have a continuous singular part. We propose a  criterion under the assumption that the continuous singular part of $F$ is dominated by the discrete part in a certain sense. The criterion generalizes the previous results by Alexeev and Khartov for discrete probability laws and the results by Berger and Kutlu for the mixtures of discrete and absolutely continuous laws. In addition, we describe the characteristic triplet of the corresponding L\'evy-type representation, which may contain some continuous singular part. We also show that the assumption of the dominated continuous singular part cannot 	be simply omitted or even slightly extended (without some special assumptions). We apply the general criterion to some interesting particular  examples. We also positively solve the decomposition problem stated by Lindner, Pan and  Sato  within the considered case.
\end{abstract}

\textit{Keywords and phrases}: distribution functions, characteristic functions, continuous singular part, infinite divisibility, rational-infinite divisibility, quasi-infinite divisibility, the L\'evy-type representation.	
	
\section{Introduction}

	This paper is devoted to the study of a new class of probability laws that naturally extends the fundamental class of infinitely divisible distributions. 
	
	Let $F$ be a  distribution function on the real line $\R$. Recall that $F$ (and the corresponding probability law) is called \textit{infinitely divisible} if for every positive integer $n$ there exists a  distribution function $F_{1/n}$ such that $F=(F_{1/n})^{*n}$, where ``$*$'' denotes the convolution operation, i.e. $F$ is the $n$-fold convolution power of $F_{1/n}$. Let $f$ be the characteristic function of $F$, i.e.
	\begin{eqnarray*}
		f(t):=\int_{\R} e^{itx} \dd F(x),\quad t\in\R.
	\end{eqnarray*}
	It is well-known (see \cite{GnedKolm}, \cite{Sato1999}, and \cite{Zolot2}) that $F$ is infinitely divisible if and only if $f$ admits \textit{the L\'evy representation}:
	\begin{eqnarray}\label{repr_f}
		f(t)=\exp\biggl\{it \gamma-\dfrac{\sigma^2t^2}{2}+\int_{\R\setminus\{0\}} \bigl(e^{itx} -1 -it \sin(x)\bigr)\dd L(x)\biggr\},\quad t\in\R,
	\end{eqnarray}
		with some \textit{shift parameter} $\gamma\in\R$, \textit{the Gaussian variance} $\sigma^2 \geqslant 0$, and \textit{the L\'evy spectral function} $L:\R\setminus\{0\}\to \R$, which is  non-decreasing on every interval $(-\infty,0)$ and $(0,+\infty)$,  and it satisfies
	\begin{eqnarray}\label{eq_Lambda_zerooninfty}
		\lim_{x \to -\infty}L(x) = \lim_{x\to+\infty}L(x)  = 0,
	\end{eqnarray}
	and also
	\begin{eqnarray*}
		\int_{O_\delta} x^2 \dd L(x)<\infty\quad\text{for any}\quad\delta>0,
	\end{eqnarray*}
	where $O_\delta:= (-\delta,0)\cup(0,\delta)$.  The function $L$ is assumed to be right-continuous at every point of the real line. It is important to note that \textit{the characteristic triplet} $(\gamma,\sigma^2, L)$ is uniquely determined by $f$ and hence by $F$. Due to this representation, the class of infinitely divisible probability laws found a lot of applications through the L\'evy processes (see \cite{Sato1999})  in the stochastic calculus   (see \cite{Appl}),  the teletraffic models (see \cite{Lifshits}), and the  insurance mathematics (see \cite{Schoutens}).
		                                                                                                                                                 	
	Let $\ID$ denote the class of all infinitely divisible distribution functions on the real line. This  class is naturally extended by the following  way. We call a distribution function $F$ (and the corresponding probability law) \textit{rational-infinitely divisible} if there exist some infinitely divisible distribution functions $F_1$ and $F_2$ such that $F_1=F*F_2$. In terms of characteristic functions, this definition is equivalent to  the formula $f(t)=f_1(t)/f_2(t)$, $t\in\R$, for the characteristic function $f$ of $F$, where $f_1$ and $f_2$ are the characteristic functions of some infinitely divisible distribution functions $F_1$ and $F_2$, respectively. We denote by $\RID$ the class of all rational-infinitely divisible distribution functions.	Since $F_2$ may be chosen as degenerate at some point $a$ here (i.e. $f_2(t)=e^{ita}$, $t\in\R$),  it is clear that, indeed, $\ID\subset \RID$. Moreover, from the difinition, it is seen that the characteristic function $f$ of any $F\in\RID$ admits \textit{the L\'evy-type representation}. Namely, if $F_1$ and $F_2$ have the characteristic triplets $(\gamma_1, \sigma_1^2,L_1)$ and $(\gamma_2, \sigma_2^2, L_2)$, respectively, then formula \eqref{repr_f} holds  with  \textit{the shift parameter} $\gamma=\gamma_1-\gamma_2\in\R$, \textit{the Gaussian variance} $\sigma^2=\sigma_1^2-\sigma_2^2$, and \textit{the spectral function} $L=L_1-L_2$. In that case,  $L$ has a bounded total variation on $\R\setminus O_\delta$ for every $\delta>0$, and, in general, it is non-monototic on the intervals $(-\infty,0)$ and $(0,+\infty)$. The function $L$ also inherits from $L_1$ and $L_2$ the right-continuity on $\R$ and  property \eqref{eq_Lambda_zerooninfty}. Moreover, 
	\begin{eqnarray*}
		\int_{O_\delta} x^2 \dd |L|(x)<\infty\quad\text{for any}\quad\delta>0,
	\end{eqnarray*}
	where we integrate over the variation of $L$. We now suppose that, conversely, $f$ admits representation \eqref{repr_f} with some $\gamma\in\R$, $\sigma^2\geqslant 0$, and  $L$ satisfying the above conditions. Following Lindner and Sato \cite{LindPanSato}, the corresponding distribution function $F$ (and the corresponding probability law)  is called \textit{quasi-infinitely divisible}. Let us fix any real $\gamma_1$ and $\gamma_2$ such that $\gamma=\gamma_1-\gamma_2$. Let us fix any non-negative $\sigma_1^2$ and $\sigma_2^2$ such that $\sigma_1^2\geqslant\sigma_2^2$ and $\sigma^2=\sigma_1^2-\sigma^2_2$. Due to the Hahn--Jordan decomposition, it is not difficult to show that there exist some canonical L\'evy spectral  functions $L_1$ and $L_2$ satisfying the usual conditions with monotonicity such that $L=L_1-L_2$. Then $f(t)=f_1(t)/f_2(t)$, $t\in\R$, where $f_1$ and $f_2$ are  represented by  canonical L\'evy's formula \eqref{repr_f} with the characteristic triplets $(\gamma_1, \sigma_1^2, L_1)$ and $(\gamma_2, \sigma_2^2, L_2)$, respectively. So $f_1$ and $f_2$ are the characteristic functions for  some infinitely divisible distribution functions and hence the distribution function $F$ corresponded to such an $f$ is rational-infinitely divisible. Thus $F\in\RID$ if and only if $f$ admits representation \eqref{repr_f} with some $\gamma$, $\sigma^2$, and $L$ satisfying the above conditions. Additionally, \textit{the characteristic triplet} $(\gamma, \sigma^2, L)$ is uniquely determined by $f$ and hence by $F$ as for infinitely divisible laws (it can be concluded from the assertion from  \cite{GnedKolm}, p.~80).  It is also clear that for any rational-infinitely divisible $F$ its characteristic function $f$ has no zeroes on the real line, i.e. $f(t)\ne 0$, $t\in\R$.
 	
	The class $\RID$ and its multivariate analog  are actively studied now (see \cite{AlexeevKhartov_Mult}, \cite{BergKutLind}, and \cite{LindPanSato}) and they find some interesting applications in probability limit and compactness theorems (see Sections~4 and~8 in \cite{LindPanSato}, Section~3 in \cite{AlexeevKhartov}, the paper \cite{KhartovWeak}, and also \cite{AlexeevKhartov2}, \cite{KhartovComp}), and in other areas (see, for instance,  \cite{ChhDemniMou}, \cite{Nakamura}, and \cite{Passegg}). But, actually, non-degenerate representatives of $\RID\setminus \ID$ appeared even earlier in the theory of decompositions of probability laws as components of some infinitely divisible distribution functions (see \cite{GnedKolm} p. 81--83 and \cite{LinOstr} p. 165).
	
	The class $\RID$ is seen to be rather wide. For instance, it contains the distribution function of every probability law, which has a mass $>1/2$ at some point. Hence the class $\RID$  contains non-degenerate distribution functions of some probability laws with  bounded supports (see examples in \cite{LindPanSato}), which are ``far''  from the   infinite divisibility property in a known sense (see \cite{BaxterShapiro}). So it is rather interesting and important to obtain criteria for belonging to the class $\RID$. The existing results in this direction usually have  quite simple and nice formulations in terms of characteristic functions. 	The first quite general result of such type was obtained by Lindner, Pan, and Sato in \cite{LindPanSato} (see  Theorem~8.1, p.~30). It states that a lattice distribution function $F$ belongs to the class $\RID$  if and only if its  characteristic function $f$ does not have  zeroes on the real line, i.e. $f(t)\ne 0$ for any $t\in\R$.
	In the articles \cite{AlexeevKhartov} and \cite{Khartov}, this result was generalized to the class of arbitarary discrete probability laws. Namely, a discrete distribution function $F$ belongs to $\RID$ if and only if its characteristic function $f$ is separeted from zero, i.e.  $|f(t)|\geqslant \mu$  for any $t\in\R$ and for some constant $\mu>0$. Moreover, in that case,  the components of the characteristic triplet  are fully described. This result is  a generalization of the previous  one for discrete lattice distributions, because the absolute value $|f(\,\cdot\,)|$ of the characteristic function $f$ of a discrete lattice distribution is a periodic continuous function on $\R$. Therefore such $f$ is zero-free on the period segment (and hence on $\R$) if any only if it is separated from zero.
	
	 In the paper \cite{Berger}, Berger considered mixtures of a degenerate law (with a non-zero coefficient) and absolutely continuous distributions. According to his result, a distribution function $F$ of such type belongs to $\RID$ if and only if   $f(t)\ne 0$, $t\in\R$. Moreover, the result describes the structure of the components of the characteristic triplet in that case. The author also formulated  more general  criterion for the case, when   the degenerate law from the previous one is  replaced by a discrete lattice distribution with characteristic function, which has no zeroes on the real line. At present, however, the most general criterion (among those that use assumptions about the type of distribution) is the following result for the mixtures of discrete and absolutely continuous probability laws, which was obtained by Berger and Kutlu in the paper \cite{BergerKutlu}. Let us formulate it with more details here. Namely, assume that $F(x)=c_d F_d(x)+c_a F_a(x)$, $x\in\R$, where $F_d$ is a discrete distribution function, $F_a$ is an absolutely continuous distribution function, $c_d>0$, $c_a\geqslant 0$, and $c_d+c_a=1$. We write the characteristic function $f$ in the corresponding form:  $f(t)=c_d f_d(t)+c_a f_a(t)$, $t\in\R$. Then $F\in\RID$ if and only if $f(t)\ne 0$ and $|f_d(t)|\geqslant \mu$ for any $t\in\R$ with some constant $\mu>0$.  It is equivalent to the condition $|f(t)|\geqslant \mu'$ for any $t\in\R$ and for some constant $\mu'>0$. Moreover,  Berger and Kutlu showed the existence of some discrete part in the spectral function and they fully described its absolutely continuous part for this case. It should be noted that we are not aware any similar results for pure absolutely continuous distribution functions $F$. However, for some cases the problem of membership to the class $\RID$ for a given  distribution function of such type is not difficult to solve by the general criteria proposed in \cite{KhartovGeneral} with some additional analysis.
	
	This article is devoted to generalization and complementation all the mentioned results (except \cite{KhartovGeneral}) for the case, when $F$ may have a continuous singular part. Namely, we propose a  criterion of belonging a distribution function $F$ to the class $\RID$ under the assumption that its continuous singular part is dominated by its discrete part in a certain sense. In fact, we show that the conditions on $f$ from the results \cite{Berger} and \cite{BergerKutlu} are carried over to this case. 	In addition, we describe the characteristic triplet of the corresponding L\'evy-type representation, which may contain some continuous singular part. 
	We next show that the assumption of the dominated continuous singular part cannot 	be simply omitted or even slightly extended without some additional assumptions. In addition, for any  $F\in\RID$ we solve the decomposition problem, which was stated by Lindner, Pan, Sato in \cite{LindPanSato} (see Open Question~8.4),  within the considered case. Here we obtain a positive solution generalizing similar rusults from \cite{Berger} and \cite{BergerKutlu}.
	
	The article has the following structure. Section~2 contains necessary preliminaries, more detailed statements of some existing results mentioned above and the formulations of the new  results of the paper.  In Section~3, we formulate some important known theorems and useful  lemmata, which will be auxiliary tools  needed for the proofs of our results. In Section~4, we first prove a key auxiliary lemma and   we next propose the proofs of the main results of the article.              
	
	Throughout the paper, we use the following notation. We denote by $\N$  the set of positive integers and let $\N_0:=\N\cup \{0\}$. The symbols $\Z$ and $\Q$ denote as usual the sets of all integers and rational numbers, respectively. Next, $\CC$ is the set of all complex numbers.  For any $z\in\CC$ we denote by $\ImagPart\{z\}$ and $\arg(z)$ its imaginary part and the principal value of the argument of $z$, respectively. If $\psi$ is a complex-valued  continuous function on $\R$ satisfying $\psi(0)=c\in\CC\setminus\{0\}$ and $\psi(t)\ne 0$ for any $t\in\R$, then  \textit{the distinguished logarithm} $\Ln\psi$ is defined by the formula $\Ln\psi(t):=\ln |\psi(t)|+i \Arg \psi(t)$, $t\in\R$, where  $\Arg \psi(t)$ is the argument of $\psi(t)$ uniquely defined on $\R$ by the continuity with the condition $\Arg \psi(0)=\arg (c)\in (-\pi,\pi]$. The symbol $\id_a$ with fixed $a\in\R$ denotes the distribution function of the degenerate law concentrated at the point $a$, i.e. $\id_a(x)=1$ for $x\geqslant a$ and $\id_a(x)=0$ for $x<a$. For any set $A$ we denote by $\idd_A$ the indicator function of $A$, i.e. $\idd_A(x)=1$ for any $x\in A$ and $\idd_A(x)=0$ for any $x\notin A$.  The signum function is denoted by $\sgn(\cdot)$, i.e. $\sgn(x)=+1$ for $x>0$, $\sgn(x)=-1$ for $x<0$, and $\sgn(0)=0$. For any finite set $A$ the symbol $|A|$ denotes the number of elements of $A$.
	We always set  $\sum_{k\in A} a_k= 0$ and $\prod_{k\in A} a_k= 1$ in the case $A=\varnothing$.
	For any two vectors $x$ and $y$ from $\R^n$ the standard scalar product is denoted by $\lng x,y\rng$.	
	
	For any function $G$ defined on $\R$ the limits $\lim_{x \to\pm\infty} G(x)$ are denoted by $G(\pm \infty)$, respectively, if these limits exist.  The class of all functions $G:\R\to\R$ of bounded total  variation on $\R$ (non-monototic in general), which are right-continuous at every point and $G(-\infty)=0$, is denoted by $\BV$. We denote by $\BVC$ the set of all functions $G:\R\to\CC$ of the form $G(x)=G_1(x)+iG_2(x)$, $x\in\R$, where $G_1$ and $G_2$ are from $\BV$.	For every $G\in\BV$ (or $\BVC$) its total variation on $\R$ will be denoted by $\|G\|$ and the total variation on $(-\infty,x]$ --- by $|G|(x)$, $x\in\R$. So we have $|G(x)|\leqslant |G|(x)\leqslant\|G\|$, $x\in\R$, and $|G|(+\infty)=\|G\|$. Next, we adopt the following convention. Let $G$ be a function from $\BV$ with the Fourier--Stieltjes transform $g$, i.e. $g(t)=\int_{\R} e^{itx}\dd G(x)$, $t\in\R$. In view of the uniqueness theorem  for functions from $\BV$, we set $\|g\|:= \|G\|$. So $\|g\|=0$ if and only if $g(t)=0$, $t\in\R$, and  $\|c\cdot g\|=|c|\cdot\|g\|$ for any $c\in\R$.  Let $G_1$ and $G_2$ be functions from $\BV$ with the Fourier--Stieltjes transforms $g_1$ and $g_2$, correspondingly. The known inequalities   $\|G_1+G_2\|\leqslant \|G_1\|+\|G_2\|$ and $\|G_1*G_2\|\leqslant \|G_1\|\cdot\|G_2\|$ are correspondingly written as   $\|g_1+g_2\|\leqslant \|g_1\|+\|g_2\|$  and $\|g_1\cdot g_2\|\leqslant \|g_1\|\cdot \|g_2\|$. 	We recall that both $\BV$  and the corresponding  space of functions $g$  with norm  $\|\,\cdot\,\|$ will be complete normed spaces (see \cite{Gelfand} p. 165).

	\section{Criteria for belonging to class $\RID$}
	
	Let $F$ be an arbitrary distribution function on the real line. According to the Lebesgue decomposition theorem,  $F$ admits the representation:
	\begin{eqnarray}\label{repr_F_Lebesgue}
		F(x) = c_d F_{d}(x) + c_a F_{a}(x)+ c_s F_s(x), \quad x\in \R, 
	\end{eqnarray}
	where $F_d$, $F_a$, and $F_s$ are discrete, absolutely continuous and continuous singular distribution functions, respectively. Here the coefficients $c_a$, $c_a$, and $c_s$  are non-negative constants such that $c_d+c_a+c_s=1$. Let $f$ be the characteristic function of $F$. It is represented in the similar way:
	\begin{eqnarray}\label{repr_f_Lebesgue}
		f(t) = c_d f_{d}(t) + c_a f_{a}(t)+ c_s f_s(t), \quad t\in \R, 
	\end{eqnarray}
	where $f_d$, $f_a$, and $f_s$ are the characteristic functions corresponding to $F_d$, $F_a$, and  $F_s$, respectively. It is well known that the summands in \eqref{repr_F_Lebesgue} or \eqref{repr_f_Lebesgue} are uniquely determined. So if any of the terms is not identically zero, then the corresponding coefficient, distribution function, and characteristic function are uniquely determined.
	
	We will consider only the case, when  $F$ has a  non-zero discrete part, i.e. $c_d>0$ in \eqref{repr_F_Lebesgue}. We write the distribution function $F_d$ in the form
	\begin{eqnarray}\label{def_Fd}
		F_{d}(x)= \sum_{\substack{k\in \N_0:\\ x_k\leqslant x}} p_k, \quad x\in \R,
	\end{eqnarray}
	where  $x_k$ are distinct reals associated with weights $p_k \geqslant 0$, $k \in \N_0$, $\sum_{k=0}^{\infty} p_k = 1$. Hence $f_d$ has the form
	\begin{eqnarray}\label{def_fd}
		f_{d}(t)= \sum_{k\in\N_0} p_k e^{i t x_k},\quad t\in \R.
	\end{eqnarray}
	 We define the support of the distribution corresponding to $F_d$:
	 \begin{eqnarray*}
	 	\X := \{ x_k:  p_{k} > 0, k \in \N_0\}.
	 \end{eqnarray*}
	 Obviously, $\X\ne \varnothing$. We also need the set of all finite $\Z$--linear combinations of elements from the set $\X$:
	 \begin{eqnarray*}
	 	\lng \X\rng:=\biggl\{\sum_{k=1}^m a_k z_k:\,   a_k \in \Z,\, z_k\in \X,\, m\in\N\biggr\}.
	 \end{eqnarray*}
	 So $\lng \X\rng$ is the module over the ring $\Z$ with the generating set $\X$. It easily seen that, in particular, $\X\subset\lng \X\rng$ and $0\in \lng \X\rng$. If $\X\ne\{0\}$, then $\lng \X\rng$ is an infinite countable set.
	 
	 Now we are ready to formulate in details the most general existing results on criteria for belonging to the class $\RID$.	 We start with the result obtained in \cite{AlexeevKhartov} and \cite{Khartov} by Alexeev and Khartov for the case of discrete $F$.
	
	\begin{Theorem}\label{th_Repr_Disc}
		Suppose that $F$ is a discrete distribution function, $c_d=1$ and $c_a=c_s=0$ in \eqref{repr_F_Lebesgue} $(${}$F$ and $F_d$ are identical and hence $f$ and $f_d$ are too;  $F$ has form \eqref{def_Fd}, $f$ is represented by \eqref{def_fd}$)$. Then $F\in\RID$  if and only if $\inf_{t\in\R} |f(t)| > 0$. In that case, $f$ admits the following representation
		\begin{eqnarray}\label{discr_Levy}
			f(t) = \exp\biggl\{ i t\gamma_0 + \sum_{u \in \langle \X \rangle \setminus\{0\}} \lambda_u \bigr(e^{i tu} - 1\bigr) \biggr\}, \quad t\in \R,
		\end{eqnarray} 
		with some $\gamma_0 \in \lng \X \rng $ and $\lambda_u \in \R$ for all $u \in \lng \X \rng \setminus\{0\}$, and $\sum_{u \in \langle \X \rangle \setminus\{0\}}|\lambda_u| < \infty$.
	\end{Theorem} 
	
	It is not difficult to rewrite  representation \eqref{discr_Levy} in  integral form \eqref{repr_f} with some $\gamma\in\R$, $\sigma^2=0$, and some discrete $L$, which will  satisfy all conditions for a spectral function for the quasi-infinitely divisibility (we will do it below for more general case). Hence, in particular, if the characteristic function $f$ is represented by \eqref{discr_Levy}, then $F\in\RID$.
	
	It is clear that if $F\in\RID$ and $\lambda_u\geqslant 0$ for all $u \in \lng \X \rng \setminus\{0\}$ in representation \eqref{discr_Levy}, then $F$ is infinitely divisible, i.e. $F\in\ID$. If there exists $\lambda_{v}<0$ with some $v\in \lng \X \rng \setminus\{0\}$, then $F\in\RID\setminus\ID$, because the (uniquely defined) function $L$ will be decreasing in the neighborhood of $v$. There are particular examples of the latter case in \cite{LindPanSato} (p. 10) and \cite{LinOstr} (p. 165).

	We now formulate the result by Berger and Kutlu from the paper \cite{BergerKutlu} for the important case $c_a\geqslant 0$. As we said before, at present, it is the most general criterion now using information about the type of the distribution function $F$.
	\begin{Theorem}\label{th_Repr_DiscAbs}
		Suppose that $F$ has decomposition \eqref{repr_F_Lebesgue} with some $c_d>0$,  $c_a\geqslant 0$, and $c_s=0$. Then the following statements are equivalent:
		\begin{enumerate}
			\item[$(i)$] $F\in\RID$,
			\item[$(ii)$] $\inf_{t\in\R} |f(t)| > 0$,
			\item[$(iii)$] $f(t)\ne 0$ for any $t\in\R$, and $\inf_{t\in\R} |f_d(t)| > 0$.
		\end{enumerate}
		If one of the conditions is satisfied, and hence all, then $f$ admits the following representation
		\begin{eqnarray}\label{discr_abs_Levy}
				f(t)=\exp\Biggl\{it\gamma_0+\sum_{u \in \mathcal{U}} \lambda_u \bigr(e^{i tu} - 1\bigr)+\int_{\R} (e^{itx}-1)  \biggl( v(x)+\sgn(x)\,  \dfrac{\mathfrak{m}\cdot e^{-|x|}}{|x|}\biggr) \dd x\Biggr\},\quad t\in\R,
		\end{eqnarray}
	where $\gamma_0\in\R$, $\mathcal{U}$ is a discrete set, $\lambda_u \in \R$ for all $u \in \mathcal{U}$, and $\sum_{u \in \mathcal{U}}|\lambda_u|<\infty$, $\mathfrak{m}\in\Z$,  $v$ is a real-valued function from $L_1(\R)$.
	\end{Theorem}
		
	It should be noted that, in addition to this theorem, the paper \cite{BergerKutlu} contains a number of interesting assertions concerning  decomposition of $F$ and  characterization of the existence of the $H$--moment of $F$ for certain positive functions $H$.  Additionally, there is a rather general theorem similar to Theorem~\ref{th_Repr_DiscAbs} for the function $f$, but without the assumption that it is a characteristic function of some probability law (Theorem~2.1 in \cite{BergerKutlu}).
		
	Theorem~\ref{th_Repr_DiscAbs} generalizes above mentioned (in introduction) results by Berger from \cite{Berger} (Theorems~4.5 and 4.12), where $F_d$ was assumed to be discrete lattice.  It is also seen that, in fact, Theorem~\ref{th_Repr_DiscAbs} is a strengthening of Theorem~\ref{th_Repr_Disc}. However, we note that $\gamma_0$ and the discrete part in \eqref{discr_Levy} are described  in a greater detail than in formula \eqref{discr_abs_Levy}.
	
	We now propose a generalization of the previous criteria for the case, when $F$ may have a continuous singular part, i.e. $c_s\geqslant 0$ in \eqref{repr_F_Lebesgue}. The following theorem is the main result of this article.
	
	For convenience, we preliminarily select the following property of distributions. 	Let $\mu_d:= \inf_{t \in \R} |f_d(t)|$. We say that a distribution function $F$ has \textit{the dominated continuous singular part} if 
	$c_s< c_d \mu_d$ for the case $\mu_d>0$ and if $c_s=0$ for the case $\mu_d=0$.

	\begin{Theorem}\label{th_MainResult}
		Suppose that $F$ has decomposition \eqref{repr_F_Lebesgue} with some $c_d>0$, $c_a\geqslant 0$, $c_s\geqslant 0$, and $F$ has the dominated continuous singular part. Then the following statements are equivalent:
		\begin{enumerate}
			\item[$(i)$] $F\in\RID$,
			\item[$(ii)$] $\inf_{t\in\R} |f(t)| > 0$,
			\item[$(iii)$] $f(t)\ne 0$ for any $t\in\R$, and $\inf_{t\in\R} |f_d(t)| > 0$.
		\end{enumerate}		
		If one of the conditions is satisfied, and hence all, then $f$ admits the following representation
		\begin{multline}\label{th_MainResult_Levy}
			f(t)=\exp\Biggl\{ i t\gamma_0  + \sum_{u \in \langle \X \rangle \setminus\{0\}} \lambda_u \bigr(e^{i tu} - 1\bigr)+\\
			+\int_{\R\setminus\{0\}} (e^{itx}-1)  \biggl( v_a(x)+\sgn(x)\,  \dfrac{\mathfrak{m}_a\cdot e^{-|x|}}{|x|}\biggr) \dd x\\
			+\int_{\R\setminus\{0\}} (e^{itx}-1)\,\dd W(x)\Biggr\},\quad t\in\R.
		\end{multline} 
		Here $\gamma_0 \in \langle \X \rangle$, $\lambda_u \in \R$ for all $u \in \lng \X \rng \setminus\{0\}$, and $\sum_{u \in \lng \X \rng \setminus\{0\}}|\lambda_u| < \infty$. Next, the function $v_a:\R\mapsto \R$ satisfies $\int_{\R} |v_a(x)|\dd x<\infty$, and, in the case $c_a=0$,  $v_a$ is identically $0$; the constant $\ind_a\in\Z$ is defined by the formula
		\begin{eqnarray}\label{def_inda}
			\ind_a := \dfrac{1}{2\pi} \Bigl(\lim_{t\to\infty} \Arg R_{a}(t) - \lim_{t\to-\infty} \Arg R_{a}(t)\Bigr),
		\end{eqnarray}
		with
		\begin{eqnarray*}
			R_a(t):=1+\dfrac{c_af_a(t)}{c_df_d(t)+c_sf_s(t)},\quad t\in\R,
		\end{eqnarray*}
		where, in particular, $\ind_a=0$ for the case $c_a=0$. Next, the function $W: \R\to\R$ belongs to $\BV$ and it is always continuous on $\R$. If $c_s=0$ then $W$ is identically $0$. If $c_s\ne 0$ then $W$ is not absolutely continuous on $\R$, i.e. it always contains some continuous singular part. In addition, if  all the functions $F_s^{*k}$, $k\in\N$, are continuous singular, then the function $W$ is $($pure\,$)$ continuous singular.
	\end{Theorem}
	It is important to note that the discrete part in the exponent in  \eqref{th_MainResult_Levy} depends only on discrete part of $f$, i.e. on $f_d$. More precisely, we have the following remark, which will be seen from the proof of Theorem~\ref{th_MainResult}.
	\begin{Remark}\label{rem_Discrete} In representation \eqref{th_MainResult_Levy}, 
		 \begin{eqnarray*}
		 	\Ln f_d(t) = i t\gamma_0 + \sum_{u \in \langle \X \rangle \setminus\{0\}} \lambda_u \bigr(e^{i tu} - 1\bigr), \quad t\in \R.
		 \end{eqnarray*} 
	\end{Remark}

	It is easily seen that if we set $c_a:=0$ and $c_s:=0$ in the assumptions of this theorem, then we come to Theorem~\ref{th_Repr_Disc}.  If we set only $c_s:=0$, then we conclude the statement of Theorem~\ref{th_Repr_DiscAbs}, but with a full discription of the discrete part  in \eqref{discr_abs_Levy}. In the unexplored case $c_s>0$, the new function $W$ appears and it always has some continuous singular part.

	Using Theorem~\ref{th_MainResult}, it is easy to construct a lot of particular examples of $F\in\RID$ with non-zero continuous singular parts.\\
	
	\noindent\textbf{Example 1.} Suppose that
	\begin{eqnarray*}
		F(x):=c_d\id_0(x)+ c_a F_a(x)+c_s F_s(x),\quad x\in\R,
	\end{eqnarray*}
	with $c_d>c_s>0$, $c_a\geqslant 0$, and $c_d+c_a+c_s=1$. Let  $F_s$ be an arbitrary continuous singular function, but $F_a$ be an absolutely continuous distribution function, whose characteristic function $f_a$ is  real and non-negative (for instance, $f_a$ is a P\'olya-type characteristic function or $f_a$ is corresponded to a symmetric continuous stable distribution). Then $F\in\RID$.  Indeed, $F$ has the dominated continuous singular part, because $f_d$ is identically $1$ here, i.e. $\mu_d=1$, and hence $c_d\mu_d=c_d>c_s$. Let us check $(ii)$ from Theorem~\ref{th_MainResult}. For this, we consider the characteristic function of $F$ 
	\begin{eqnarray*}
		f(t)=c_d+c_af_a(t)+c_s f_s(t),\quad t\in\R.
	\end{eqnarray*}
	Under the assumptions, we observe that
	\begin{eqnarray*}
		|f(t)|\geqslant|c_d+c_af_a(t)|-c_s |f_s(t)|=c_d+c_af_a(t)-c_s |f_s(t)| ,\quad t\in\R.
	\end{eqnarray*}
	Since $f_a(t)\geqslant 0$ and $|f_s(t)|\leqslant 1$ for any $t\in\R$, we have $|f(t)|\geqslant c_d -c_s>0$ for any $t\in\R$, i.e. $(ii)$  holds. Thus $F\in\RID$ by Theorem~\ref{th_MainResult}.\quad $\Box$\\

	In general case, of course, checking  $(ii)$ or $(iii)$ of Theorem~\ref{th_MainResult} may require  more subtle analisys.  We illustrate it by the following rather special example, which is also interesting, because there is a (pure) continuous singular function $W$ from \eqref{th_MainResult_Levy}.\\
	
	\noindent\textbf{Example 2.} Let us consider the mixture of the degenerate law, the uniform distribution on $[0,1]$ and the classical Cantor distribution. Namely, we set
	\begin{eqnarray*}
		F(x):=c_d\id_0(x)+ c_a U(x)+c_s S(x),\quad x\in\R,
	\end{eqnarray*}
	where $U$ is an absolutely continuous distribution function with the density $\idd_{[0,1]}$, and $S$ is the cumulative function of the classical Cantor distribution supported on the Cantor set  $\mathcal{C}\subset[0,1]$; $c_d>c_s>0$, $c_a\geqslant 0$, and $c_d+c_a+c_s=1$.  
	
	We first observe that $f_d$ is identically $1$ here and, similarly as in Example~1, $F$ has the dominated continuous singular part. We next check condition $(iii)$ of Theorem~\ref{th_MainResult}.  Since $\inf_{t\in\R} |f_d(t)|=1$, it remains to  show that $f(t)\ne 0$ for any $t\in\R$. The function $f$ is expressed by the formula
	\begin{eqnarray*}
		f(t)=c_d+c_a\,\dfrac{e^{it}-1}{it}+ c_s e^{it/2} \prod_{k=1}^\infty \cos (t/3^k),\quad t\in\R.
	\end{eqnarray*}
	We set
	\begin{eqnarray*}
		\hat{f}(t):= f(t) e^{-it/2}=c_d e^{-it/2}+c_a\,\dfrac{\sin(t/2)}{t/2}+ c_s \prod_{k=1}^\infty \cos (t/3^k),\quad t\in\R.
	\end{eqnarray*}
	The functions $f$ and $\hat{f}$ have the same set of zeroes on $\R$. Observe that it contains in the set of zeroes of the function $\ImagPart\{\hat{f}(t)\}=-c_d \sin(t/2)$, $t\in\R$. The last set is exactly $\{2\pi m: m\in\Z\}$. Therefore it is sufficient to show that $f(2\pi m)\ne 0$ for any $m\in\Z$. Obviously, we can exclude $m=0$. For any $m\in\Z\setminus\{0\}$ we write
	\begin{eqnarray*}
		|f(2\pi m)|=|\hat{f}(2\pi m)|&=&\biggl| c_d \cos(\pi m)+ c_a\,\dfrac{\sin(\pi m)}{\pi m}+ c_s \prod_{k=1}^\infty \cos (2\pi m/3^k) \biggr|\\
		&=&\biggl| c_d (-1)^m+  c_s \prod_{k=1}^\infty \cos (2\pi m/3^k) \biggr|.
	\end{eqnarray*}
	It follows that
		\begin{eqnarray*}
		|f(2\pi m)|\geqslant  c_d- c_s \prod_{k=1}^\infty \bigl|\cos (2\pi m/3^k) \bigr|\geqslant c_d-c_s>0.
	\end{eqnarray*}
		Thus $f(t)\ne 0$ for any $t\in\R$ and condition $(iii)$ of Theorem~\ref{th_MainResult} is satisfied. So, by this theorem, $F\in\RID$.
	
	Let us consider representation \eqref{th_MainResult_Levy} for $f$. Since $f_d(t)=1$ for any $t\in\R$,   the sum $i t\gamma_0 + \sum_u  \lambda_u \bigr(e^{i tu} - 1\bigr)$ is identically $0$  according to Remark~\ref{rem_Discrete}. There are the function $v_a\in L_1(\R)$ and the constant $\mathfrak{m}_a\in\Z$, but we will not determine these here. Next, we turn to the function $W$. It is known that all convolution powers $S^{*n}$, $n\in\N$, of the Cantor distribution function $S$ are continuous singular. We recall that this fact follows from the famous Jessen--Wintner theorem (see \cite{JessenWintner}, Theorem 35) and the representation of $S(x+\tfrac{1}{2})$, $x\in\R$, as the infinite symmetric Bernoulli convolution, but it was explicitly shown in  \cite{WienerWintner} (see p.~520). According to Theorem~\ref{th_MainResult}, it implies that $W$ is pure continuous singular. Moreover, it will be seen from the proof of Theorem~\ref{th_MainResult} that $W$ is the limit  of the sums
	\begin{eqnarray*}
		W_n(x):=\sum_{k=1}^{n}  \dfrac{(-1)^{k-1} }{k}\biggl(\dfrac{c_s }{c_d}\biggr)^k S^{*k}(x),\quad x\in\R,\quad n\in\N,
	\end{eqnarray*} 
	as $n\to\infty$ with respect to convergence in total variation on $\R$.\quad $\Box$\\

	So Theorem~\ref{th_MainResult} yields the sufficient condition for $W$ to be (pure) continuous singular. Is it true that $W$ is actually always (pure) continuous singular function in the case $c_s\ne 0$?
	
	\begin{Proposition} \label{pr_PureSing}
		Suppose that $F\in\RID$ satisfies the assumptions of Theorem~\ref{th_MainResult} with $c_s>0$, and  the characteristic function $f$ is represented by \eqref{th_MainResult_Levy} with some $W$. Suppose that there exists an integer $n_a\geqslant 2$ such that the function $F_s^{*(n_a-1)}$ is $($pure\,$)$ continuous singular, but the function $F_s^{*n_a}$ is not, i.e.  $F_s^{*n_a}(x)=\alpha H_a(x)+(1-\alpha)H_s(x)$, $x\in\R$, where $\alpha$ is a number from $(0,1]$, $H_a$ is an absolutely continuous distribution function, $H_s$ is a continuous singular distribution function. If $\alpha\geqslant \tfrac{n_a}{n_a+1}\cdot\tfrac{c_s}{c_d}$ then the function $W$ is not $($pure\,$)$ continuous singular. In particular, it is always true for the case $\alpha =1$.
	\end{Proposition}
	
	It should be recalled here that there exist examples of continuous singular distribution functions (say of $F_s$), whose convolution squares (say $F_s^{*2}$) are absolutely continuous (see \cite{WienerWintner} and also  \cite{HewittZuckermann}, \cite{IvashevMusatov}). Therefore we answer the question asked before Proposition~\ref{pr_PureSing} in the negative.

	The following remark yields the characteristic triplet for $F\in\RID$ in the explicit form under the assumptions of Theorem~\ref{th_MainResult}.
	\begin{Remark}
		Representation \eqref{th_MainResult_Levy} can be written in the following  forms:
		\begin{eqnarray*}
			f(t)=\exp\biggl\{ i t\gamma_0  + \int_{\R\setminus\{0\}} (e^{itx}-1)\,\dd L(x)\biggr\}=\exp\biggl\{ i t\gamma  + \int_{\R\setminus\{0\}} \bigl(e^{itx}-1-it\sin(x)\bigr)\,\dd L(x)\biggr\}
			,\quad t\in\R,
		\end{eqnarray*}
		where 
		\begin{eqnarray*}
			\gamma:=\gamma_0+ \int_{\R\setminus\{0\}} \sin(x)\,\dd L(x)
		\end{eqnarray*}
		and $L(x):=L_d(x)+L_a(x)+L_s(x)$ for any $x\in\R\setminus\{0\}$ with 
	\begin{eqnarray}
		L_d(x)&:=&
		\begin{cases}
			\sum\limits_{\substack{u \in \langle \X \rangle \setminus\{0\}:\\
					u\leqslant x}} \lambda_u,& x<0,\\
			-\sum\limits_{\substack{u \in \langle \X \rangle \setminus\{0\}:\\
					u> x}}  \lambda_u,& x>0,
		\end{cases}\label{def_Ld}\\
		L_a(x)&:=&
		\begin{cases}
			\int\limits_{u\leqslant x} \Bigl( v_a(u)+\sgn(u)\,  \tfrac{\mathfrak{m}_a\cdot e^{-|u|}}{|u|}\Bigr) \dd u,& x<0,\\
			-\int\limits_{u> x}  \Bigl( v_a(u)+\sgn(u)\,  \tfrac{\mathfrak{m}_a\cdot e^{-|u|}}{|u|}\Bigr) \dd u,& x>0,
		\end{cases}\label{def_La}\\
		L_s(x)&:=&
		\begin{cases}
			W(x),& x<0,\\
			W(x)-W(+\infty),& x>0.
		\end{cases}	\label{def_Ls}
	\end{eqnarray}
	\end{Remark}
	
It is seen that $L$ satisfies all admissible conditions for a spectral function in the L\'evy type representation. Thus, under the assumptions of Theorem~\ref{th_MainResult}, $(\gamma,0,L)$ is the characteristic triplet for $F$ satisfying one of the conditions $(i)$--$(iii)$. On the other hand, if we know that \eqref{th_MainResult_Levy} represents the characteristic function of some probability law, then its distribution function $F$ is quasi-infinitely divisible by the definition and hence $F\in\RID$. 
	
The following remark shows that the property of the dominated singular part is not necessary condition for belonging to the class $\RID$.

\begin{Remark}\label{rem_notdomin}
	 Let $F$ satisfy the assumptions of Theorem~\ref{th_MainResult} with $c_s>0$. Suppose that $F\in\RID$.  Then $F^{*n}\in\RID$ for any $n\in\N$, but  $F^{*n}$ don't have the dominated singular parts for all sufficiently large $n$.
\end{Remark}
Indeed, it is seen from the definition of the rational-infinite divisibility that the convolution of  two any distribution functions from $\RID$ belongs to $\RID$. Therefore $F^{*n}\in\RID$ for any $n\in\N$. Next, we consider the characteristic function $f$ of $F$. It admits decomposition \eqref{repr_f_Lebesgue} with $c_d>0$, $c_a\geqslant 0$, $c_s>0$, and $c_s<c_d\mu_d$, where $\mu_d=\inf_{t \in \R}|f_d(t)|$ as before. Then every $F^{*n}$ has the following characteristic function:
\begin{eqnarray*}
	f(t)^n=\bigl(c_df_d(t)+c_af_a(t)+c_sf_s(t) \bigr)^n,\quad t\in\R.
\end{eqnarray*}
From this it is easily seen that $t\mapsto c_d^n f_d(t)^n$, $t\in\R$, is the  Fourier--Stieltjes transform of the discrete part of $F^{*n}$. Hence $c_d^n$ is the weight of this part and $\inf_{t \in \R}|f_d(t)^n|=\mu_d^n$. Next, there are $n$ terms $c_sc_d^{n-1} f_s(t)f_d(t)^{n-1}$, which are included to the  Fourier--Stieltjes transform of the continuous singular part of $F^{*n}$. This is due to the well known fact that a convolution of continuous singular function with a discrete function is continuous singular (see \cite{Lukacs} p. 190 or \cite{Tucker} p. 319). So the weight of the continuous singular part of $F^{*n}$ is not less than $n c_s c_d^{n-1}$. For any integer $n\geqslant c_d/c_s$ we have that $n c_s c_d^{n-1}\geqslant c_d^n\geqslant c_d^n\mu_d^n$. Therefore the exact weight of the continuous singular part of $F^{*n}$ is not less than $c_d^n\mu_d^n$ too, i.e. the condition of the dominated singular part  doesn't hold for $F^{*n}$ for $n\geqslant c_d/c_s$.

The following interesting example shows that the condition of the dominated singular part can not be simply omitted and, moreover, it cannot be extended to the case $c_s=c_d\mu_d$ with $\mu_d>0$ without certain additional assumptions.\\

\noindent\textbf{Example 3.} Let $B$ denote the distribution function of the Bernoulli law on the points $\pm 1$ with equal probabilities, i.e. $B(x)=\tfrac{1}{2}\,\id_{-1}(x)+\tfrac{1}{2}\,\id_{1}(x)$, $x\in\R$. We set 
\begin{eqnarray}\label{def_F*}
	F_*(x):=\dfrac{1}{2}\,\id_0(x)+\dfrac{1}{2}\, F_s(x),\quad x\in\R,
\end{eqnarray}
where $F_s$ is the following infinite symmetric Bernoulli convolution:
\begin{eqnarray*}
	F_s=B_1*B_2*\ldots*B_n*\ldots
\end{eqnarray*}
with $B_k(x):=B(k!\, x)$, $x\in\R$, $k\in\N$. Let $f_*$ and $f_s$ denote the characteristic functions of $F_*$ and $F_s$, respectively. Then $f_*(t)=\tfrac{1}{2}+\tfrac{1}{2}\, f_s(t)$, $t\in\R$, and
\begin{eqnarray}\label{def_fs}
	f_s(t)=\prod\limits_{k=1}^{\infty} \cos (t/k!),\quad  t\in\R.
\end{eqnarray} 
It is known (see \cite{Lukacs} pp. 20 and 67) that the function $F_s$ is continuous singular. Let us consider  Lebesgue's decomposition \eqref{repr_f_Lebesgue} for $f=f_*$. We have $c_d=\tfrac{1}{2}$, $c_a=0$, and $c_s=\tfrac{1}{2}$. Observe that the component $f_d$ is identically $1$ for $f_*$ and hence $\mu_d=\inf_{t \in \R}|f_d(t)|=1$. So $c_d\mu_d= c_s$, i.e. the condition of the dominated continuous singular part doesn't hold. Next, $f_*(t)\ne 0$ for any $t\in\R$, since otherwise $f_s(t_0)=-1$ for some $t_0\ne 0$ and hence $|f_s(t_0)|=1$, which would mean that $F_s$ is a discrete lattice distribution function (see  \cite{Lukacs} Theorem 2.1.4), a contradiction. Thus $F_*$ satisfies condition $(iii)$ of Theorem~\ref{th_MainResult} and the assumption $c_d>0$, but $F_*$ doesn't have the dominated continuous singular part.    
     
\begin{Proposition}\label{pr_Example}
	The function $F_*$ doesn't belong to $\RID$.
\end{Proposition}
Thus, in general,  conditions $(i)$ and $(iii)$ from Theorem~\ref{th_MainResult} are not equivalent even if $c_s=c_d\mu_d>0$.\quad $\Box$\\

By the way, we recall that distribution functions of discrete laws with a point mass $\tfrac{1}{2}$, which have characteristic functions without zeroes on the real line, don't always belong to the class $\RID$ too (see \cite{KhartovAlexeevThree} for more details). However, if the distribution is an atom of mass $\tfrac{1}{2}$ plus a continuous part, which is not purely singular, then the answer is definite.\\

\noindent\textbf{Example 4.} Suppose that
\begin{eqnarray*}
	F(x):= \dfrac{1}{2}\,\id_{\gamma_0}(x)+ c_a F_a(x)+c_s F_s(x),\quad x\in\R,
\end{eqnarray*}
where $F_a$ is an absolutely continuous distribution function, $F_s$ is a continuous singular distribution function, $\gamma_0\in\R$, $c_a>0$, $c_s\geqslant 0$, and  $c_a+c_s=\tfrac{1}{2}$.  Such a function $F$ always belongs to the class $\RID$. Let us check it. We consider the characteristic function of $F$:
\begin{eqnarray*}
	f(t)=\dfrac{1}{2}\, e^{it\gamma_0} +c_a f_a(t)+c_s f_s(t),\quad t\in\R.
\end{eqnarray*}
Here $c_d=\tfrac{1}{2}$ and $f_d(t)=e^{it\gamma_0}$, $t\in\R$. Hence $\mu_d=\inf_{t\in\R}|f_d(t)|=1$ and  $c_s=\tfrac{1}{2}-c_a<\tfrac{1}{2}=c_d\mu_d$, i.e. $F$ has the dominated continuous singular part.  Next, $f(t)\ne 0$ for any $t\in\R$. Indeed,   otherwise $2c_af_a(t_0)+2c_s f_s(t_0)=-e^{it_0\gamma_0}$ for some $t_0\ne 0$ and, in particular, $|2c_af_a(t_0)+2c_s f_s(t_0)|=1$. Then the latter implies that the distribution function with the characteristic function $2c_af_a(\cdot)+2c_s f_s(\cdot)$ is  discrete lattice, but this is actually continuous by the assumptions, a contradiction.

So $F$ satisfies condition $(iii)$ and other assumptions of Theorems~\ref{th_MainResult}. It follows that $F\in\RID$. \quad $\Box$\\

The next assertion solves the decomposition problem for any  distribution function $F\in\RID$ satisfying the assumptions of Theorem~\ref{th_MainResult}.	
\begin{Proposition}\label{pr_Decom}
	Let $F$ be a distribution function with decomposition \eqref{repr_F_Lebesgue} with some $c_d>0$, $c_a\geqslant 0$, $c_s\geqslant 0$, and the dominated continuous singular part. Suppose that $F\in\RID$ and $F=F_1*F_2$ with some distribution functions $F_1$ and $F_2$.  Then $F_1$ and $F_2$ belong to $\RID$.
\end{Proposition}
	
	This proposition generalizes the similar result from \cite{BergerKutlu} (see Corollary 2.3 from Theorem 2.2), which solves the problem for $F\in\RID$ satisfying the assumptions of Theorem~\ref{th_Repr_DiscAbs}. We note that there is a general decomposition problem for arbitrary $F\in\RID$, which was stated by Lindner, Pan, Sato in \cite{LindPanSato} (see Open Question~8.4): \textit{is it true that if $F\in\RID$ and $F=F_1*F_2$ $(${}$F_1$ and $F_2$ are distribution functions on $\R${}$)$, then $F_1\in\RID$ and $F_2\in\RID$?} So Proposition~\ref{pr_Decom} answers this question in the positive for any  $F\in\RID$ satisfying the assumptions of Theorem~\ref{th_MainResult}.	However, there is a result in \cite{KhartovDecomQ}, which asserts that the general answer is negative.

	\section{Auxiliary theorems and lemmata}
	
	The proof of the main result (Theorem~\ref{th_MainResult}) of the article essentially uses the following Wiener--Pitt theorem (see \cite{WienerPitt}, \cite{Shreider}, and \cite{Gelfand}, p. 191).
	
	\begin{Theorem}\label{th_WienerPitt}
		Let $H$ be a function from $\BVC$ with the following Lebesgue's decomposition
		\begin{eqnarray*}
			H(x)=H_d(x)+H_a(x)+H_s(x),\quad x\in\R,
		\end{eqnarray*}
		where $H_d$, $H_a$, $H_s$ are correspondingly  discrete, absolutely continuous and singular parts of $H$, which belong to $\BVC$. Let
		\begin{eqnarray*}
			h(t):=\int_{\R} e^{itx} \dd H(x),\quad\text{and}\quad h_d(t):= \int_{\R} e^{itx} \dd H_d(x),\quad t\in\R.
		\end{eqnarray*}
	Suppose that $\inf_{t \in \R}|h(t)|>0$ and $\|H_s\|< \inf_{t \in \R}|h_d(t)|$. Then there exists a function $K$ from $\BVC$ such that
	\begin{eqnarray*}
		\dfrac{1}{h(t)}=\int_{\R} e^{itx} \dd K(x),\quad t\in\R.
	\end{eqnarray*}
	\end{Theorem}
	
	The following theorem was proposed by Berger \cite{Berger}. In fact, it is a modification of one Krein's result (see \cite{Krein}, Theorem L, p.15).
		
	\begin{Theorem}\label{th_Berger}
		Suppose that a function $h: \R\to\CC$ is defined by the formula
		\begin{eqnarray}\label{def_h_Berger}
			h(t):= c+ \int_{\R}e^{itx} u(x)\dd x,\quad t\in\R,
		\end{eqnarray}
		where $c\in\CC\setminus\{0\}$ is a constant, $u:\R\to\CC$ is a function satisfying the condition $\int_{\R} |u(x)|\dd x<\infty$. Assume that $h(t)\ne 0$ for any $t\in\R$. Then
		\begin{eqnarray*}
			h(t)=\exp\Biggl\{\int_{\R\setminus\{0\}} (e^{itx}-1)  \biggl( v(x)+\sgn(x)\,  \dfrac{\mathfrak{m}\cdot e^{-|x|}}{|x|}\biggr) \dd x\Biggr\},\quad t\in\R,
		\end{eqnarray*}
	where $\ind$ is a constant, which is defined by the formula  
	\begin{eqnarray*}
		\ind := \dfrac{1}{2\pi} \Bigl(\lim_{t\to\infty} \Arg h(t) - \lim_{t\to-\infty} \Arg h(t)\Bigr),
	\end{eqnarray*}
	and the function $v:\R\to\CC$ satisfies the condition $\int_{\R} |v(x)|\dd x<\infty$.
	\end{Theorem}
	
	Let us consider the quantity $\ind$ from the theorem. It is the \textit{index} of $h$ (see \cite{Krein}). Observe that, by the Riemann--Lebesgue lemma, for the function $h$ defined by \eqref{def_h_Berger} we have $h(t)\to c$ as $t\to\pm\infty$. Therefore it is not difficult to conclude the following remark.
	\begin{Remark}\label{rem_ind}
		The quantity  $\ind$ from Theorem~\ref{th_Berger} is well-defined and it is an integer.
	\end{Remark}

	We next formulate a very useful assertion obtained by Berger \cite{Berger}.
	
	\begin{Theorem}\label{th_Berger2}
	Let $F$ be a distribution function on $\R$ with the characteristic function $f$. Suppose that $f$ admits the following representation
	\begin{eqnarray*}
		f(t)=\exp\biggl\{it \gamma_1-\dfrac{\sigma^2t^2}{2}+\int_{\R\setminus\{0\}} \Bigl(e^{itx} -1 -it\,\idd_{[-1,1]}(x)\Bigr)\dd L(x)\biggr\},\quad t\in\R,
	\end{eqnarray*}
	with some  $\gamma_1\in\CC$,  $\sigma^2 \in\CC$, and the function $L(x)=L_1(x)+iL_2(x)$, $x\in\R\setminus\{0\}$. Here for every $j\in\{1,2\}$ the function $L_j$ has a bounded variation on $\R\setminus O_\delta$ for any $\delta>0$ $($in general, $L_j$ may be non-monototic on intervals $(-\infty,0)$ and $(0,+\infty)${}$)$, $L_j$ is right-continuous at every point on $\R$, $L_j(-\infty)=L_j(+\infty)=0$, and  
	\begin{eqnarray*}
		\int_{O_\delta} x^2 \dd |L_j|(x)<\infty\quad\text{for any}\quad\delta>0,
	\end{eqnarray*}
	where $O_\delta:= (-\delta,0)\cup(0,\delta)$. Then $\gamma_1\in\R$, $\sigma^2\geqslant 0$, and $\ImagPart\{L(x)\}=L_2(x)=0$ for any $x\in\R\setminus\{0\}$. In addition, $F\in\RID$.
	\end{Theorem}
	
	The next theorem yields one special property of characteristic functions of rational-infinitely divisible laws. Its proof can be found in  \cite{Khartov} (p. 3) or \cite{KhartovWeak} (p. 360).
	\begin{Theorem}\label{th_psifunction}
		Let $F$ be a distribution function on $\R$ with the characteristic function $f$. If $F\in\RID$ then for any $\tau>0$ there exists  $C_\tau>0$ such that 
		\[
		\sup_{t\in\R}\biggl|\dfrac{f(t-\tau)f(t+\tau)}{f(t)^2}\biggr|\leqslant C_\tau.
		\]
	\end{Theorem}

	The following lemma can be found in the paper by Bochner \cite{Bochner}.
	
	\begin{Lemma}\label{lm_Bochner}
		Let $(W_n)_{n\in\N}$ be a sequence from $\BV$ such that $\sup_{n\in\N}\|W_n\|=B<\infty$. Let
		\begin{eqnarray*}
			w_n(t):=\int_{\R} e^{itx} \dd W_n(x), \quad t\in \R, \quad n\in\N.
		\end{eqnarray*}
		Suppose that $w_n(t)\to w(t)$ as $n\to\infty$ uniformly on every bounded interval.  Then there exists $W\in \BV$ with $\|W\|\leqslant B$ such that
		\begin{eqnarray*}
			w(t)=\int_{\R} e^{itx} \dd W(x),\quad t\in \R.
		\end{eqnarray*}
	\end{Lemma}

	We will need the following multivariate version of the Leibniz rule for the partial derivatives of a product of two functions (see \cite{Hardy} p. 10). 
	\begin{Lemma}\label{lm_partialderiv_productrule}
		Suppose that functions $f:\R^d\to \CC$ and $g:\R^d\to \CC$  have continuous mixed partial derivatives up to and including order $d$  in some points $t=(t_1,\ldots, t_d)$  from some open set $T\subset\R^d$.  Then for any $t\in T$
		\begin{eqnarray*}
			\dfrac{\partial^d (f(t) \cdot g(t)) }{\partial t_1\ldots \partial t_d}=\sum_{\subs \subset D}\Biggl( \dfrac{\partial^{|\subs| } f(t)}{\prod_{j\in\subs} \partial t_j}\cdot  \dfrac{\partial^{d-|\subs| }g(t)}{\prod_{j\in D\setminus\subs} \partial t_j}\Biggr),
		\end{eqnarray*}
		where $D=\{1,\ldots, d\}$, and the parameter $\subs$ of the sum assumes all subsets of $D$. 
	\end{Lemma}
	
	We will also need  the multivariate form of Fa\`a di Bruno's formula for the partial derivatives of a composition of two functions (see \cite{Hardy} p. 4).
	\begin{Lemma}\label{lm_partialderiv_composition}
		Suppose that a function $g:\R^d\to \CC$  has continuous mixed partial derivatives up to and including order $d$ in some points $t=(t_1,\ldots, t_d)$ from some open set $T\subset\R^d$. Suppose that $f:\CC\to \CC$ is a holomorphic function in some domain in $\CC$, which include the values of $g(t)$, $t\in T$.  Suppose that $D:=\{1,\ldots,d\}$, $\subs$ is non-empty subset of $D$, $\mathcal{P}(\subs)$ is the set of all partititions of the set $\subs$. Then for any $t\in T$
		\begin{eqnarray*}
			\dfrac{\partial^{|\subs|} f (g(t)) }{\prod_{j\in\subs}\partial t_j}=\sum_{P\in \mathcal{P}(\subs)} \Biggl(\dfrac{\dd^{|P|} f(z)}{\dd z^{|P|}}\Bigg|_{z= g(t)}\cdot \prod_{s \in P} \dfrac{\partial^{|s| } g(t)}{\prod_{j\in s} \partial t_j}\Biggr).
		\end{eqnarray*}
	\end{Lemma}

	We now formulate a lemma, which will be a key tool in proving of Theorem~\ref{th_MainResult}.  The next section begins with the proof of this lemma.	
	\begin{Lemma}\label{lm_summodulkoefZ}
		Suppose that a function $\varphi: \R^d\to \CC$, where $d\in\N$, admits the following representation
		\begin{eqnarray*}
			\varphi(t)=\sum_{m=1}^{N} q_{m} e^{i\lng t, c_m\rng},\quad t=(t_1,\ldots,t_d)\in\R^d,  
		\end{eqnarray*}
		with $N\in\N$, $q_m\in \CC$, and $c_m=(c_{m,1},\ldots, c_{m,d})\in\Z^d$, $m=1,\ldots, N$. Then for any $k\in\N$
			\begin{eqnarray}\label{lm_summodulkoefZ_ineq}
			\|\varphi^k\|\leqslant  \biggl(\dfrac{\pi}{\sqrt{6}}\biggr)^d   \sum_{\subs \subset D} \Biggl[\,\prod_{j\in D\setminus\subs} |k\alpha_j-1|\cdot \sum_{\substack{P\in \mathcal{P}(\subs):\\ |P|\leqslant k}} \biggl(\, \prod_{l=1}^{|P|} (k-l+1)\cdot S_\varphi^{k-|P|} \cdot R_\varphi(P)\biggr) \Biggr],
		\end{eqnarray}
	where $D=\{1,\ldots, d\}$,
	\begin{eqnarray}\label{lm_summodulkoefZ_defSR}
		&&\qquad\qquad\alpha_j:=\min_{m=1,\ldots, N} c_{m,j},\quad j=1,\ldots, d,\nonumber\\
		&&S_\varphi:= \sup\limits_{t\in[-\pi,\pi]^d}|\varphi(t)|,\qquad R_\varphi(P):=\sup\limits_{t\in[-\pi,\pi]^d}\prod_{s \in P} \biggl|\dfrac{\partial^{|s| }\varphi(t)}{\prod_{j\in s} \partial t_j}\biggr|.
	\end{eqnarray}
	In particular, 
	\begin{eqnarray}\label{lm_summodulkoefZ_ineq2}
		\|\varphi^k\|\leqslant A_\varphi k^d S_\varphi^k\quad\text{for any}\quad k\in\N
	\end{eqnarray}
	with
	\begin{eqnarray}\label{lm_summodulkoefZ_defAphi}
		A_\varphi:=\biggl(\dfrac{\pi}{\sqrt{6}}\biggr)^d  
		\sum_{\subs \subset D} \Biggl[\,\prod_{j\in D\setminus\subs} \bigl(|\alpha_j|+1\bigr)
		\cdot\!\!\! \sum_{P\in \mathcal{P}(\subs)} S_\varphi^{-|P|}R_\varphi(P)\Biggr],
	\end{eqnarray}
	which is independent of $k$.
	\end{Lemma}

\section{Proofs}
\textbf{Proof of Lemma \ref{lm_summodulkoefZ}.} We will first prove inequality \eqref{lm_summodulkoefZ_ineq} for $k=1$. We define a vector	 $a:=(a_1,\ldots, a_d)$ with $a_j:= 1-\alpha_j$, $j=1,\ldots, d$,	and we introduce the function
\begin{eqnarray*}
	\tphi_a(t):=\varphi(t)e^{i\lng t,a\rng}=\sum_{m=1}^{N} q_{m} e^{i\lng t, c_m+a\rng},\quad t\in\R^d.		
\end{eqnarray*}
Note that $c_m+a\in\N^d$ for any $m=1,\ldots,N$. Without loss of generality we can assume that  $c_m$ are distinct if $N>1$. Then we have
\begin{eqnarray}
	\|\varphi\|=\sum_{m=1}^N |q_m|&=& \sum_{m=1}^N \biggl(|q_m|\cdot \prod\limits_{j=1}^{d}(c_{m,j}+a_j)\cdot\prod\limits_{j=1}^{d}\dfrac{1}{c_{m,j}+a_j} \biggr)\nonumber\\
	&\leqslant&\biggl(\,\sum_{m=1}^N \biggl[|q_m|^2\cdot\prod\limits_{j=1}^{d}(c_{m,j}+a_j)^2\biggr] \biggr)^{1/2}\cdot\biggl(\,\sum_{m=1}^N  \prod\limits_{j=1}^{d}\dfrac{1}{(c_{m,j}+a_j)^2}\biggr)^{1/2}.\label{ineq_norm_phi}
\end{eqnarray}
Since $c_m+a=(c_{m,1}+a_1,\ldots, c_{m,d}+a_d)\in\N^d$ are distinct vectors for different $m$, the following estimate holds	
\begin{eqnarray*}
	\sum_{m=1}^N  \prod\limits_{j=1}^{d}\dfrac{1}{(c_{m,j}+a_j)^2}\leqslant \prod\limits_{j=1}^{d}\biggl(1+\dfrac{1}{2^2}+\ldots+\dfrac{1}{n^2}+\ldots\biggr)=\biggl(\dfrac{\pi^2}{6}\biggr)^d.
\end{eqnarray*}	
We next write the function $\tphi_a$ in the expanded form:
\begin{eqnarray*}
	\tphi_a(t)=\sum_{m=1}^{N} q_{m}\exp\Bigl\{i\sum_{j=1}^{d}t_j (c_{m,j}+a_j)\Bigr\},\quad t=(t_1,\ldots,t_d)\in\R^d.
\end{eqnarray*}
So it is easily seen that
\begin{eqnarray*}
	\dfrac{\partial^d \tphi_a (t) }{\partial t_1\ldots \partial t_d}= \sum_{m=1}^{N} \biggl(q_{m}\cdot i^d\cdot\prod\limits_{j=1}^{d}(c_{m,j}+a_j)\cdot\exp\Bigl\{i\sum_{j=1}^{d}t_j (c_{m,j}+a_j)\Bigr\}\biggr),\quad t=(t_1,\ldots,t_d)\in\R^d.
\end{eqnarray*}
We consider the written  exponential functions (additionally normed by $(2\pi)^d$) as $N$-rank orthonormal system of $L_2\bigl([-\pi,\pi]^d\bigr)$ space and, by Parseval's identity, we get
\begin{eqnarray*}
	\sum_{m=1}^N \biggl[|q_m|^2\cdot\prod\limits_{j=1}^{d}(c_{m,j}+a_j)^2\biggr]=\dfrac{1}{(2\pi)^d} \int_{[-\pi,\pi]^d} \biggl|\dfrac{\partial^d \tphi_a (t) }{\partial t_1\ldots \partial t_d} \biggr|^2 \dd t_1\ldots\dd t_d.
\end{eqnarray*}
Hence
\begin{eqnarray*}
	\sum_{m=1}^N \biggl[|q_m|^2\cdot\prod\limits_{j=1}^{d}(c_{m,j}+a_j)^2\biggr]\leqslant \sup\limits_{t\in[-\pi,\pi]^d} \biggl|\dfrac{\partial^d \tphi_a (t) }{\partial t_1\ldots \partial t_d} \biggr|^2 .
\end{eqnarray*}
Using the obtained estimates in \eqref{ineq_norm_phi}, we come to the inequality
\begin{eqnarray*}
	\|\varphi\|\leqslant  \biggl(\dfrac{\pi}{\sqrt{6}}\biggr)^d \sup\limits_{t\in[-\pi,\pi]^d} \biggl|\dfrac{\partial^d \tphi_a (t) }{\partial t_1\ldots \partial t_d} \biggr|.
\end{eqnarray*}
We next find an upper estimate for the latter supremum in terms of $\varphi$ and $a$.    We apply Lemma~\ref{lm_partialderiv_productrule} for the function $\tphi_a(t)=\varphi(t)e^{i\lng t,a\rng}$ with $f(t):=\varphi(t)$ and $g(t):= e^{i\lng t,a\rng}$, $t\in\R^d$:  
\begin{eqnarray*}
	\dfrac{\partial^d \tphi_a (t) }{\partial t_1\ldots \partial t_d}=\sum_{\subs \subset D}\Biggl( \dfrac{\partial^{|\subs| }\varphi(t)}{\prod_{j\in\subs} \partial t_j}\cdot  \dfrac{\partial^{d-|\subs| }e^{i\lng t, a\rng}}{\prod_{j\in D\setminus\subs} \partial t_j}\Biggr)=\sum_{\subs \subset D}\Biggl( \dfrac{\partial^{|\subs| }\varphi(t)}{\prod_{j\in\subs} \partial t_j}\cdot e^{i\lng t, a\rng} \cdot i^{d-|\subs|} \prod_{j\in D\setminus\subs}  a_j \Biggr),\quad t\in\R^d.
\end{eqnarray*}
Here $D:= \{1,\ldots, d\}$ and the index $\subs$ of the sum  assumes all subsets of $D$ as values. Hence we get the estimate:
\begin{eqnarray*}
	\sup\limits_{t\in[-\pi,\pi]^d} \biggl|\dfrac{\partial^d \tphi_a (t) }{\partial t_1\ldots \partial t_d} \biggr|\leqslant \sum_{\subs \subset D} \Biggl(  \prod_{j\in D\setminus\subs} | a_j|\cdot\sup\limits_{t\in[-\pi,\pi]^d}\Biggl| \dfrac{\partial^{|\subs| }\varphi(t)}{\prod_{j\in\subs} \partial t_j}\Biggr| \Biggr).
\end{eqnarray*}
Thus (on account of equalities $a_j=1-\alpha_j$, $j=1,\ldots, d$)  we have
\begin{eqnarray}\label{ineq_norm_phi1}
	\|\varphi\|\leqslant  \biggl(\dfrac{\pi}{\sqrt{6}}\biggr)^d  \sum_{\subs \subset D} \Biggl(  \prod_{j\in D\setminus\subs} |\alpha_j-1|\cdot\sup\limits_{t\in[-\pi,\pi]^d}\Biggl| \dfrac{\partial^{|\subs| }\varphi(t)}{\prod_{j\in\subs} \partial t_j}\Biggr| \Biggr).
\end{eqnarray}
Observe that the right-hand side of this inequality coincides with the right-hand side of \eqref{lm_summodulkoefZ_ineq} for $k=1$. Indeed, in the inner sum in \eqref{lm_summodulkoefZ_ineq}, the index $P$ can be equal to $\{\subs\}$, i.e. $|P|=1$, and, in the product in $R_\varphi(P)$, the index $s$ assumes only the value  $\subs$. Therefore $\prod_{l=1}^{|P|} (k-l+1)= 1$, $S_\varphi^{k-|P|}=S_\varphi^0=1$, and   
\begin{eqnarray*}
	R_\varphi(P)= \sup\limits_{t\in[-\pi,\pi]^d}\Biggl| \dfrac{\partial^{|\subs| }\varphi(t)}{\prod_{j\in\subs} \partial t_j}\Biggr|.
\end{eqnarray*}
Substituting these in \eqref{lm_summodulkoefZ_ineq}, we see the matching with \eqref{ineq_norm_phi1}.

We now prove that \eqref{lm_summodulkoefZ_ineq} is true for any $k\in\N$.	Let us fix arbitrary $k\in\N$. We consider the function $\varphi^k$:
\begin{eqnarray*}
	\varphi^k(t)=\biggl(\sum_{m=1}^{N} q_{m} e^{i\lng t, c_m\rng}\biggr)^k,\quad t\in\R^d.
\end{eqnarray*}
Clearly, it can be written in the form
\begin{eqnarray*}
	\varphi^k(t)= \sum_{m=1}^{M} Q_m e^{i\lng t,C_m\rng},\quad t\in\R^d,
\end{eqnarray*}
with some $M\in\N$, $Q_m\in \CC$, and distinct $C_m=(C_{m,1},\ldots, C_{m,d})\in\Z^d$, $m=1,\ldots, M$. Observe that
\begin{eqnarray*}
	C_m\in\Bigl\{\sum_{j=1}^k c_{m_j}: m_1,\ldots, m_k\in \{1,\ldots, N\} \Bigr\},\quad m=1,\ldots, M.
\end{eqnarray*}
Therefore, in particular,
\begin{eqnarray*}
	\min_{m=1,\ldots, M} C_{m,j}= k \min_{m=1,\ldots, N} c_{m,j}=k\alpha_j,\quad j=1,\ldots, d.
\end{eqnarray*} 
Taking this into account, we use inequality \eqref{ineq_norm_phi1} for the function $\varphi^k$:
\begin{eqnarray*}
	\|\varphi^k\|\leqslant\biggl(\dfrac{\pi}{\sqrt{6}}\biggr)^d  \sum_{\subs \subset D} \Biggl(\prod_{j\in D\setminus\subs} |k\alpha_j-1|\cdot\sup\limits_{t\in[-\pi,\pi]^d}\Biggl| \dfrac{\partial^{|\subs| }\varphi^k(t)}{\prod_{j\in\subs} \partial t_j}\Biggr|   \Biggr).
\end{eqnarray*}
We now write expressions for the mixed partial derivatives from the right-hand side.   We apply Lemma~\ref{lm_partialderiv_composition} with $g(t):= \varphi(t)$, $t\in\R^d$, and  $f(z):=z^k$, $z\in\CC$:
\begin{eqnarray*}
	\dfrac{\partial^{|\subs| }\varphi^k(t)}{\prod_{j\in\subs} \partial t_j}
	&=&\sum_{P\in \mathcal{P}(\subs)} \Biggl(\dfrac{\dd^{|P|} z^k}{\dd z^{|P|}}\Bigg|_{z=\varphi(t)}\cdot \prod_{s \in P} \dfrac{\partial^{|s| }\varphi(t)}{\prod_{j\in s} \partial t_j}\biggr)\\
	&=&\sum_{\substack{P\in \mathcal{P}(\subs):\\ |P|\leqslant k}} \Biggl( \prod_{l=1}^{|P|}(k-l+1)\cdot \varphi(t)^{k-|P|} \cdot \prod_{s \in P} \dfrac{\partial^{|s| }\varphi(t)}{\prod_{j\in s} \partial t_j}\Biggr).
\end{eqnarray*}
We note that this formula holds including the case $\subs=\varnothing$, where $P$ must be $\varnothing$, i.e. $|P|=0$ and $\mathcal{P}(\subs)=\{\varnothing\}$. Indeed, there is only one summand in the latter sum, in which the written products are equal to $1$ (because the indexes formally belong the empty sets).  Therefore the right-hand side is equal to $\varphi(t)^k$  as the left-hand side since $|\subs|=0$.  

We have the estimate
\begin{eqnarray*}
	\sup\limits_{t\in[-\pi,\pi]^d}\Biggl| \dfrac{\partial^{|\subs| }\varphi^k(t)}{\prod_{j\in\subs} \partial t_j}\Biggr|\leqslant \sum_{\substack{P\in \mathcal{P}(\subs):\\ |P|\leqslant k}} \biggl( \prod_{l=1}^{|P|}(k-l+1)\cdot S_\varphi^{k-|P|} \cdot R_\varphi(P)\biggr),
\end{eqnarray*}
where $S_\varphi$ and $R_\varphi(P)$ are defined by formulas \eqref{lm_summodulkoefZ_defSR}.
Thus we obtain the required inequality
\begin{eqnarray}\label{lm_summodulkoefZ_ineq'}
	\|\varphi^k\|\leqslant \biggl(\dfrac{\pi}{\sqrt{6}}\biggr)^d  
	\sum_{\subs \subset D} \Biggl[\, \prod_{j\in D\setminus\subs} | k\alpha_j-1|
	\cdot \sum_{\substack{P\in \mathcal{P}(\subs):\\ |P|\leqslant k}} \biggl(\,\prod_{l=1}^{|P|}(k-l+1)\cdot S_\varphi^{k-|P|} \cdot R_\varphi(P)\biggr)\Biggr].
\end{eqnarray}

We next show that \eqref{lm_summodulkoefZ_ineq2} holds.  Let us obtain an upper estimate for the right-hand side of inequality \eqref{lm_summodulkoefZ_ineq'}. We fix $k\geqslant 1$. Observe that
\begin{eqnarray*}
	\prod_{j\in D\setminus\subs} | k\alpha_j-1|\leqslant \prod_{j\in D\setminus\subs} \bigl(k|\alpha_j|+1\bigr)\leqslant \prod_{j\in D\setminus\subs} \bigl(k|\alpha_j|+k\bigr)\leqslant k^{d-|\subs|}\cdot \prod_{j\in D\setminus\subs} \bigl(|\alpha_j|+1\bigr).
\end{eqnarray*}
Next, since $|P|\leqslant |\subs|$ for any $P\in\mathcal{P}(\subs)$, we have
\begin{eqnarray*}
	\prod_{l=1}^{|P|}(k-l+1)\leqslant \prod_{l=1}^{|P|}(k-1+1)=k^{|P|}\leqslant k^{|\subs|}.
\end{eqnarray*}
We apply these inequalities to \eqref{lm_summodulkoefZ_ineq'}:
\begin{eqnarray*}
	\|\varphi^k\|&\leqslant& \biggl(\dfrac{\pi}{\sqrt{6}}\biggr)^d  
	\sum_{\subs \subset D} \Biggl[\,k^{d-|\subs|}\cdot \prod_{j\in D\setminus\subs} \bigl(|\alpha_j|+1\bigr)
	\cdot \sum_{\substack{P\in \mathcal{P}(\subs):\\ |P|\leqslant k}} \bigl(\,k^{|\subs|}\cdot S_\varphi^{k-|P|} \cdot R_\varphi(P)\bigr)\Biggr]\\
	&=& k^d S_\varphi^k \cdot \biggl(\dfrac{\pi}{\sqrt{6}}\biggr)^d  
	\sum_{\subs \subset D} \Biggl[\,\prod_{j\in D\setminus\subs} \bigl(|\alpha_j|+1\bigr)
	\cdot\!\!\! \sum_{\substack{P\in \mathcal{P}(\subs):\\ |P|\leqslant k}}  S_\varphi^{-|P|} R_\varphi(P)\Biggr].
\end{eqnarray*}
Since the quantities $S_\varphi$ and $R_\varphi(P)$ are always non-negative, we come to the needed inequality
\begin{eqnarray*}
	\|\varphi^k\|\leqslant k^d S_\varphi^k \cdot \biggl(\dfrac{\pi}{\sqrt{6}}\biggr)^d  
	\sum_{\subs \subset D} \Biggl[\,\prod_{j\in D\setminus\subs} \bigl(|\alpha_j|+1\bigr)
	\cdot\!\!\! \sum_{P\in \mathcal{P}(\subs)} S_\varphi^{-|P|} R_\varphi(P)\Biggr]=A_\varphi k^d S_\varphi^k.\quad \Box
\end{eqnarray*}

\textbf{Proof of Theorem \ref{th_MainResult}.} We assume  that $F$ admits  decomposition \eqref{repr_F_Lebesgue} with some $c_d\in (0,1]$,  and $F$ has the dominated singular part. We denote $\mu_d:= \inf_{t \in \R} |f_d(t)|$. 

$(i)\Rightarrow(ii)$. Suppose that $F\in\RID$. If $\mu_d=0$ then $c_s=0$ by the assumption of the dominated singular part, and hence we can apply Theorem \ref{th_Repr_DiscAbs}. According to this theorem, $(i)$ implies $(ii)$. Next, if $\mu_d>0$ then we know that $c_s<c_d\mu_d$, i.e. $c_d\mu_d-c_s>0$. Observe that for any $t\in\R$
\begin{eqnarray*}
	|f(t)|=|c_d f_d(t)+ c_s f_s(t)+c_a f_a(t) |\geqslant c_d|f_d(t)|-c_s |f_s(t)|-c_a|f_a(t)|\geqslant c_d\mu_d-c_s-|f_a(t)|.
\end{eqnarray*}
Since $f_a(t)\to 0$ as $t\to\pm \infty$, there exists $T>0$ such that $|f_a(t)|< \tfrac{1}{2}\,(c_d\mu_d -c_s)$ as $|t|>T$. Hence $|f(t)|>\tfrac{1}{2}\,(c_d\mu_d -c_s)>0$ as $|t|>T$. Next, let us consider the function $t\mapsto |f(t)|$ on the segment $[-T,T]$. Since it is continuous, there exists $t_{\min}\in[-T,T]$ such that $|f(t)|\geqslant|f(t_{\min})|$ for any $t\in[-T,T]$. Due to the condition $F\in\RID$, we know that $f(t)\ne 0$ for any $t\in\R$ (see introduction) and, in particular,  $C_T:=|f(t_{\min})|>0$. So we get $|f(t)|\geqslant C_T>0$ for any $t\in[-T,T]$. Thus  for any $t\in\R$
\begin{eqnarray*}
	|f(t)|\geqslant \min\bigl\{ \tfrac{1}{2}\,(c_d\mu_d -c_s),C_T \bigr\}>0,
\end{eqnarray*}
i.e. $\inf_{t\in\R} |f(t)|>0$. So we come to $(ii)$.

$(ii)\Rightarrow(iii)$. Obviously, $(ii)$ yields that $f(t)\ne 0$ for any $t\in\R$. To obtain a contradiction, suppose that $\mu_d=\inf_{t\in\R}|f_d(t)|=0$. Then $c_s=0$ by the condition of the dominated singular part. According to Theorem~\ref{th_Repr_DiscAbs}, $(iii)$ follows from $(ii)$, i.e., in particular, we have that $\mu_d>0$, a contradiction.  Thus we conclude that $\inf_{t\in\R}|f_d(t)|>0$ and $(iii)$ holds.

$(iii)\Rightarrow(i)$. Let us consider $f$ represented by formula \eqref{repr_f_Lebesgue} with $c_d>0$. Here we assume that $f(t)\ne 0$ for any $t\in\R$, 
\begin{eqnarray}\label{ineq_inffd}
	\mu_d = \inf_{t \in \R}|f_d(t)|>0,
 \end{eqnarray}
and $c_s<c_d\mu_d$. Due to these assumptions, for any $t\in\R$ we have
\begin{eqnarray*}
	|c_d f_d(t)+  c_s f_s(t)|\geqslant c_d|f_d(t)|-c_s |f_s(t)|\geqslant c_d\mu_d-c_s>0,
\end{eqnarray*}
i.e.
\begin{eqnarray}\label{ineq_inffdfs}
	\inf_{t \in \R}|c_d f_d(t)+  c_s f_s(t)|>0.
\end{eqnarray}
Then for any $t\in\R$ we can write
\begin{eqnarray*}
	f(t)&=&c_d f_d(t)+ c_s f_s(t)+c_a f_a(t)\\
	&=& f_d(t)\cdot \dfrac{c_d f_d(t)+ c_s f_s(t)+c_a f_a(t)}{c_d f_d(t)+  c_s f_s(t)}\cdot \dfrac{c_d f_d(t)+  c_s f_s(t)}{f_d(t)}\\
	&=&f_d(t)\cdot \biggl(1+ \dfrac{c_a f_a(t)}{c_d f_d(t)+  c_s f_s(t)}\biggr)\cdot  \dfrac{c_d f_d(t)+  c_s f_s(t)}{f_d(t)}.
\end{eqnarray*}
So it is convenient to represent 
\begin{eqnarray}\label{repr_fd_fads_fsd}
	f(t)=f_d(t)f_{a,ds}(t)f_{s,d}(t),\quad t\in\R,
\end{eqnarray}
where
\begin{eqnarray*}
	f_{a,ds}(t):= c_d+c_s+ \dfrac{c_a(c_d+c_s) f_a(t)}{c_d f_d(t)+  c_s f_s(t)},\qquad f_{s,d}(t):= \dfrac{c_d f_d(t)+  c_s f_s(t)}{(c_d+c_s)f_d(t)},\quad t\in\R.
\end{eqnarray*}

Let us consider $f_d$ represented by \eqref{def_fd}. Due to \eqref{ineq_inffd}, by Theorem~\ref{th_Repr_Disc}, it admits the representation
\begin{eqnarray}\label{repr_fd}
	f_d(t) = \exp\biggr\{ i t\gamma_0  + \sum_{u \in \langle \X \rangle \setminus\{0\}} \lambda_u \bigr(e^{i tu} - 1\bigr) \biggr\}, \quad t\in \R,
\end{eqnarray} 
with some $\gamma_0 \in \lng \X \rng $ and $\lambda_u \in \R$ for all $u \in \lng \X \rng \setminus\{0\}$, and $\sum_{u \in \langle \X \rangle \setminus\{0\}}|\lambda_u| < \infty$.

We now consider the function $f_{a,ds}$. If $c_a=0$ then $f_{a,ds}(t)=c_d+c_s=1-c_a=1$. We next suppose that $c_a>0$. The function $t\mapsto c_d f_d(t)+  c_s f_s(t)$, $t\in\R$, is separated from $0$ according to \eqref{ineq_inffdfs}  and for its continuous
 singular part we have
\begin{eqnarray*}
	\|c_s f_s \|=c_s<c_d\mu_d=c_d \inf_{t\in\R} |f_d(t)|.
\end{eqnarray*}
By Theorem~\ref{th_WienerPitt}, there exists  a function $I_{ds}:\R\to\CC$ from $\BVC$ such that
\begin{eqnarray*}
	\dfrac{1}{c_d f_d(t)+  c_s f_s(t)}=\int_{\R} e^{itx} \dd I_{ds}(x),\quad t\in\R.
\end{eqnarray*}
We next observe that
\begin{eqnarray*}
	\dfrac{ f_a(t)}{c_d f_d(t)+  c_s f_s(t)}=\int_{\R} e^{itx} \dd F_a(x)\cdot\int_{\R} e^{itx} \dd I_{ds}(x)=\int_{\R} e^{itx} \dd (F_a*I_{ds})(x),\quad t\in\R,
\end{eqnarray*}
and we write
\begin{eqnarray*}
	\dfrac{c_a(c_d+c_s) f_a(t)}{c_d f_d(t)+  c_s f_s(t)}=\int_{\R} e^{itx} \dd \widetilde{F}_a(x),\quad t\in\R,
\end{eqnarray*}
where $\widetilde{F}_a(x):= c_a(c_d+c_s)(F_a*I_{ds})(x)$, $x\in\R$.  The function $\widetilde{F}_a\in\BVC$ inherits the absolute continuity from $F_a$, i.e. there exists a function $\widetilde{p}_a:\R\to\CC$ such that $\int_{\R} |\widetilde{p}_a(x)|\dd x<\infty$ and $\widetilde{F}_a(x)=\int_{y\leqslant x} \widetilde{p}_a(y)\dd y$, $x\in\R$. Hence
\begin{eqnarray*}
	\int_{\R} e^{itx} \dd \widetilde{F}_a(x)=\int_{\R} e^{itx}  \widetilde{p}_a(x)\dd x,\quad t\in\R,
\end{eqnarray*}
and we have
\begin{eqnarray*}
	f_{a,ds}(t)=c_d+c_s+ \dfrac{c_a(c_d+c_s) f_a(t)}{c_d f_d(t)+  c_s f_s(t)}=c_d+c_s+ \int_{\R} e^{itx}  \widetilde{p}_a(x)\dd x,\quad t\in\R,
\end{eqnarray*}
where $c_d+c_s\geqslant c_d>0$. Since $f(t)\ne 0$ for any $t\in\R$ and the function $t\mapsto c_d f_d(t)+  c_s f_s(t)$, $t\in\R$, is bounded,  we conclude that $f_{a,ds}(t)\ne 0$ for any  $t\in\R$. By Theorem~\ref{th_Berger}, $f_{a,ds}$ admits the representation
\begin{eqnarray}\label{repr_fads}
	f_{a,ds}(t)=\exp\Biggl\{\int_{\R\setminus\{0\}} (e^{itx}-1)  \biggl( v_a(x)+\sgn(x)\,  \dfrac{\mathfrak{m}_a\cdot e^{-|x|}}{|x|}\biggr) \dd x\Biggr\},\quad t\in\R,
\end{eqnarray}
where $v_a:\R\to\CC$ is a function  satisfying the condition $\int_{\R} |v_a(x)|\dd x<\infty$, $\ind_a$ is a constant  defined by the formula  
\begin{eqnarray*}
	\ind_a := \dfrac{1}{2\pi} \Bigl(\lim_{t\to\infty} \Arg f_{a,ds}(t) - \lim_{t\to-\infty} \Arg f_{a,ds}(t)\Bigr).
\end{eqnarray*}
Since $c_d+c_s$ is positive, formula \eqref{def_inda} gives an equivalent definition of $\ind_a$. By the way, it is easily seen that $\ind_a=0$ if $c_a=0$. In general, $\ind_a\in\Z$ (see Remark~\ref{rem_ind}).  Therefore, in the case $c_a=0$, we set $v_a(x):=0$ for every $x\in\R$ and the formula \eqref{repr_fads} remains valid. 

We now turn to the function
\begin{eqnarray*}
	f_{s,d}(t)=\dfrac{c_d f_d(t)+  c_s f_s(t)}{(c_d+c_s)f_d(t)}=\biggl(1+\dfrac{c_s f_s(t)}{c_d f_d(t)}\biggr)\cdot \dfrac{c_d}{c_d+c_s},\quad t\in\R.
\end{eqnarray*}
From \eqref{ineq_inffdfs} we know that $f_{s,d}(t)\ne 0$ for any $t\in\R$ and we see that $f_{s,d}(0)=1$. Hence the distinguished logarithm $\Ln f_{s,d}$ is uniquely defined on $\R$ with condition $\Ln f_{s,d}(0)=0$. Due to the assumptions, 
\begin{eqnarray}\label{ineq_fractioncdfdcsfs}
	\biggl|\dfrac{c_s f_s(t)}{c_d f_d(t)}\biggr|\leqslant \dfrac{c_s}{c_d\mu_d}<1\quad\text{for any} \quad t\in\R,
\end{eqnarray}
and we can write
\begin{eqnarray}\label{eq_Lnfsd}
	\Ln f_{s,d}(t)=\ln\biggl(1+ \dfrac{c_s f_s(t)}{c_d f_d(t)}\biggr)-\ln \biggl(1+ \dfrac{c_s}{c_d}\biggr),\quad t\in\R,
\end{eqnarray} 
where $\ln(\cdot)$ returns the principal value of the logarithm.  Let us consider the first logarithm in the right-hand side. Due to \eqref{ineq_fractioncdfdcsfs}, the following expansion holds:	
\begin{eqnarray}\label{eq_lnfsfd}
	\ln\biggl(1+ \dfrac{c_s f_s(t)}{c_d f_d(t)}\biggr)=\sum_{k=1}^{\infty} \dfrac{(-1)^{k-1} }{k}\biggl(\dfrac{c_s f_s(t)}{c_d f_d(t)}\biggr)^k=\sum_{k=1}^{\infty} \dfrac{(-1)^{k-1} }{k}\biggl(\dfrac{c_s }{c_d}\biggr)^k f_s(t)^k g_d(t)^k,\quad t\in\R,
\end{eqnarray}
where we set $g_d(t):=1/f_d(t)$, $t\in\R$. We note that the series  uniformly converges on $\R$ by the Weierstrass M-test, because, due to \eqref{ineq_fractioncdfdcsfs}, for any $t\in\R$ the absolute values of its terms  are majorized by the terms of the following convergent numerical series
\begin{eqnarray*}
	\sum_{k=1}^{\infty} \dfrac{1}{k}\biggl(\dfrac{c_s}{c_d \mu_d}\biggr)^k.
\end{eqnarray*}  

Let us consider the function $g_d$ and show that it is actually a Fourier--Stieltjes transform of some function from class $\BV$. Indeed, according to  \eqref{repr_fd},
\begin{eqnarray*}	
	g_d(t)=  e^{-it\gamma_0}\cdot e^{-\lambda_0}\cdot \exp\Bigl\{-\!\!\sum_{u \in \langle \X \rangle \setminus\{0\}} \lambda_u e^{itu} \Bigr\},\quad t\in\R, 
\end{eqnarray*}  
where we set
\begin{eqnarray*}
	\lambda_0:= -\!\!\sum_{u \in \langle \X \rangle \setminus\{0\}} \lambda_u.
\end{eqnarray*}
There is an absolutely convergent Fourier series with real coefficients in the last exponential function. Hence this function itself is expanded into absolutely convergent Fourier series with the coefficients from $\R$  and the exponents from $\langle \X \rangle$. It remains true after multiplying this series by $e^{-\lambda_0}$ and $e^{-it\gamma_0}$ (that leads to $g_d$), because $\lambda_0\in \R$ and $\gamma_0\in \langle \X \rangle$. Thus we have
\begin{eqnarray}\label{conc_gdFourier}
	g_d(t)=\sum_{y \in \langle \X \rangle} q_y e^{it y},\quad t\in\R,  
\end{eqnarray} 
where $q_y\in\R$ for every $y\in \langle \X \rangle $ and $\sum_{y \in \langle \X \rangle} |q_y|<\infty$. Obviously, the series \eqref{conc_gdFourier} can be written as a Fourier--Stieltjes transform:
\begin{eqnarray*}
	g_d(t)=\int_{\R}e^{itx} \dd I_d(x),\quad t\in\R,
\end{eqnarray*}
with the discrete function $I_d\in\BV$:
\begin{eqnarray}\label{def_Id}
	I_d(x):= \sum_{y \in \langle \X \rangle:\,  y\leqslant x} q_y,\quad x\in\R.
\end{eqnarray}

Let us return to the formula \eqref{eq_lnfsfd}. We observe that for any  $k\in\N$ the functions $f_s^k$, $g_d^k$, $f_s^k\cdot g^k$ are the Fourier--Stieltjes transforms of functions from $\BV$:
\begin{eqnarray}\label{eq_fsgd}
	f_s(t)^k=\int_{\R} e^{itx}\dd F_s^{*k}(x),\quad g_d(t)^k=\int_{\R} e^{itx} \dd I_d^{*k}(x),\quad 	f_s(t)^k g_d(t)^k=\int_{\R} e^{itx} \dd (F_s^{*k}*I_d^{*k})(x),\quad t\in\R.
\end{eqnarray}
Consequently, the partial sums of the series \eqref{eq_lnfsfd} are the Fourier--Stieltjes transforms  for the following functions:
\begin{eqnarray}\label{def_Wn}
	W_n(x):=\sum_{k=1}^{n}  \dfrac{(-1)^{k-1} }{k}\biggl(\dfrac{c_s }{c_d}\biggr)^k  \bigl(F_s^{*k}*I_d^{*k}\bigr)(x),\quad x\in\R,\quad n\in\N.
\end{eqnarray} 
Let us show that the whole sum of  the series \eqref{eq_lnfsfd} admits the representation:
\begin{eqnarray}\label{eq_sumintW}
	\sum_{k=1}^{\infty} \dfrac{(-1)^{k-1} }{k}\biggl(\dfrac{c_s }{c_d}\biggr)^k f_s(t)^k g_d(t)^k=\int_{\R} e^{itx}\dd W(x),\quad t\in\R,
\end{eqnarray}
with some function $W\in\BV$. According to Lemma~\ref{lm_Bochner}, for this it is sufficient that $1)$ the partial sums of \eqref{eq_lnfsfd} uniformly converge to the infinite sum on every bounded interval and $2)$ $\sup_{n\in\N}\|W_n\|<\infty$.  The uniform convergence of the partial sums  (even on whole real line) were showed above, i.e. $1)$  holds. Let us prove $2)$.  For any $n\in\N$ we have the estimate	
\begin{eqnarray*}
	\|W_n\|\leqslant\sum_{k=1}^{n}  \dfrac{1}{k}\biggl(\dfrac{c_s }{c_d}\biggr)^k  \bigl\|F_s^{*k}*I_d^{*k}\bigr\|\leqslant\sum_{k=1}^{n}\biggl\{ \dfrac{1}{k}\biggl(\dfrac{c_s }{c_d}\biggr)^k  \bigl\|F_s^{*k}\bigr\|\cdot \bigl\|I_d^{*k}\bigr\|\biggr\}.
\end{eqnarray*} 
Here $\bigl\|F_s^{*k}\bigr\|=1$ for every $k\in\N$, because all $F_s^{*k}$ are distribution functions. Also, it is convenient for us to represent $\|I_d^{*k}\|=\|g_d^k\|$, $k\in\N$. Thus we come to the inequalities
\begin{eqnarray}\label{ineq_normWn}
	\|W_n\|\leqslant\sum_{k=1}^{n}  \dfrac{1}{k}\biggl(\dfrac{c_s }{c_d}\biggr)^k  \bigl\|g_d^k\bigr\|,\quad n\in\N.
\end{eqnarray} 

We now estimate $\|g_d^k\|$ for every $k\in\N$. Let us return to decomposition \eqref{conc_gdFourier} for the function $g_d$. The set $\langle \X\rangle$ is countable, and $\sum_{y \in \lng \X \rng} |q_y|<\infty$. So we fix arbitrary $\e\in(0,1)$ and choose  $N_\e\in\N$ and the distinct points $y_1,\ldots, y_{N_\e}\in \lng \X \rng$ such that 
\begin{eqnarray*}
	\sum_{y \in \lng \X \rng\setminus\{y_1,\ldots, y_{N_\e}\}} |q_y|<\e.
\end{eqnarray*} 
We introduce the polynomial
\begin{eqnarray}\label{def_poly}
	\poly_{\e}(t):=\sum_{m=1}^{N_\e} q_{y_m} e^{ity_m},\quad t\in\R.
\end{eqnarray} 
So we have $\|g_d-\poly_{\e} \|<\e$. We also observe that
\begin{eqnarray}\label{ineq_norm_gd-poly}
	\biggl\|\dfrac{g_d-\poly_{\e}}{g_d}\biggr\|\leqslant \|g_d-\poly_{\e} \|\cdot \biggl\| \dfrac{1}{g_d}\biggr\|=\|g_d-\poly_{\e}\|\cdot \|f_d\|=\|g_d-\poly_{\e}\|<\e.
\end{eqnarray}
Let us fix $k\in\N$ and write
\begin{eqnarray*}
	\|g_d^k\|=\biggl\|\poly_{\e}^k \cdot \biggl(\dfrac{g_d}{\poly_{\e}}\biggr)^k\biggr\|\leqslant \|\poly_{\e}^k\|\cdot \biggl\|\dfrac{g_d}{\poly_{\e}}\biggr\|^k.
\end{eqnarray*}
In order to get a convenient representation for $g_d/\poly_{\e}$, we observe that for any $n\in\N$
\begin{eqnarray*}
	1-\biggl(\dfrac{g_d(t)-\poly_{\e}(t)}{g_d(t)}\biggr)^n&=&\biggl(1-\dfrac{g_d(t)-\poly_{\e}(t)}{g_d(t)}\biggr)\cdot \sum_{j=0}^{n-1} \biggl(\dfrac{g_d(t)-\poly_{\e}(t)}{g_d(t)}\biggr)^j\\
	&=&\dfrac{\poly_{\e}(t)}{g_d(t)}\cdot \sum_{j=0}^{n-1} \biggl(\dfrac{g_d(t)-\poly_{\e}(t)}{g_d(t)}\biggr)^j,\quad t\in\R.
\end{eqnarray*}
So we have
\begin{eqnarray*}
	\dfrac{g_d(t)}{\poly_{\e}(t)}
	=\sum_{j=0}^{n-1} \biggl(\dfrac{g_d(t)-\poly_{\e}(t)}{g(t)}\biggr)^j+\biggl(\dfrac{g_d(t)-\poly_{\e}(t)}{g_d(t)}\biggr)^n\cdot \dfrac{g_d(t)}{\poly_{\e}(t)},\quad t\in\R,\quad n\in\N.
\end{eqnarray*}
Hence, due to \eqref{ineq_norm_gd-poly}, we get
\begin{eqnarray*}
	\biggl\|\dfrac{g_d}{\poly_{\e}}\biggr\|\leqslant\sum_{j=0}^{n-1} \biggl\|\dfrac{g_d-\poly_{\e}}{g_d}\biggr\|^j+\biggl\|\dfrac{g_d-\poly_{\e}}{g_d}\biggr\|^n\cdot \biggl\|\dfrac{g_d}{\poly_{\e}}\biggr\|<\sum_{j=0}^{n-1} \e^j+\e^n\cdot \biggl\|\dfrac{g_d}{\poly_{\e}}\biggr\|,\quad n\in\N.
\end{eqnarray*}
Since $\e\in(0,1)$, letting $n\to\infty$ yields
\begin{eqnarray*}
	\biggl\|\dfrac{g_d}{\poly_{\e}}\biggr\|\leqslant\sum_{j=0}^{\infty}\e^j=\dfrac{1}{1-\e}.
\end{eqnarray*}
Thus we have
\begin{eqnarray}\label{ineq_gk}
	\|g_d^k\|\leqslant\dfrac{ \|\poly_{\e}^k\|}{(1-\e)^k}.
\end{eqnarray}
So the estimation of  $\|I_d^{*k}\|=\|g_d^k\|$ is reduced to finding  an upper bound for $\|\poly_{\e}^k\|$.  

Let us consider $\poly_{\e}$ defined by formula \eqref{def_poly}. Let $Y_\e:=\{y_1, y_2,\ldots, y_{N_\e}\}$. Suppose that $Y_\e=\{0\}$, i.e.  $\sum_{y \in \lng \X \rng\setminus\{0\}} |q_y|<\e$ and $\poly_{\e}(t)=q_0$ for any $t\in\R$. Then $\|\poly_{\e}^k\|=|q_0|^k$  for any $k\in\N$. Observe that
\begin{eqnarray*}
	|q_0|=|\poly_{\e}(0)|\leqslant |g_{d}(0)|+|\poly_{\e}(0)-g_{d}(0)|,
\end{eqnarray*}
where $|g_{d}(0)|=|1/f_d(0)|=1$ and
\begin{eqnarray*}
	|\poly_{\e}(0)-g_{d}(0)|=\biggl|  \sum_{y \in \lng \X \rng\setminus\{0\}} q_y\biggr|\leqslant \sum_{y \in \lng \X \rng\setminus\{0\}} |q_y|<\e.
\end{eqnarray*}
Thus $|q_0|<1+\e$ and hence $\|\poly_{\e}^k\|<(1+\e)^k$. 

We next assume that $Y_\e\ne \{0\}$, i.e. $Y_\e$ contains some non-zero  elements. So we select a basis over $\Q$ in the set $Y_\e$  (see \cite{Levitan} p. 67--68), i.e. we choose  non-zero elements $\beta_1, \ldots, \beta_{d} \in  Y_\e$ with some $d\in\{1,\ldots, N_\e\}$, which are linearly independent over $\Q$ and for any $m \in\{1, \ldots, N_\e\}$ there exist some $r_{m,1},\ldots, r_{m,d} \in \Q$ such that $y_m = \sum_{l= 1}^{d} r_{m,l} \beta_l$.  The linear independence over $\Q$ means that the equality $r_1 \beta_1+\ldots +r_d \beta_d=0$  holds with  $r_1,\ldots, r_d\in\Q$ only in the case  $r_1=r_2=\ldots=r_d=0$. It is clear that the coefficients $r_{m,l}$ of the decomposition of $y_m$ are uniquely determined. Let $\varkappa$ be the minimal positive integer such that $\tr_{m,l}:= \varkappa\cdot r_{m,l}\in\Z$ for any admissible  $m$ and $l$. We set $\tbeta_{l}:= \beta_l/ \varkappa$ for every $l\in\{1, \ldots, d\}$. Then we have $y_m = \sum_{l= 1}^{d} \tr_{m,l} \tbeta_l$ for any $m \in\{1, \ldots, N_\e\}$. Here the coefficients $\tr_{m,l}$ are uniquely determined by $y_m$ too. We define the vectors $\tr_m:= (\tr_{m,1},\ldots, \tr_{m,d})\in\Z^d$, $m=1,\ldots, N_\e$. Since $y_1,\ldots, y_{N_\e}$ are assumed to be distinct, all $\tr_m$ are distinct too. Let us introduce the function
\begin{eqnarray*}
	\varphi_\e(t_1,\ldots, t_d):=\sum_{m=1}^{N_\e} q_{y_m} \exp\biggl\{i\sum\limits_{l=1}^d \tr_{m,l} t_l\biggr\},\quad t_1,\ldots, t_d \in \R.
\end{eqnarray*}
It is easily seen that this function is continuous and $2\pi$-periodic over every variable. We will also use the short expression for it:
\begin{eqnarray*}
	\varphi_\e(t)=\sum_{m=1}^{N_\e} q_{y_m} e^{i\lng t,\tr_m\rng} ,\quad t=(t_1,\ldots, t_d) \in \R^d.
\end{eqnarray*}
The functions $\varphi_\e$ and $\poly_{\e}$ are related:
\begin{eqnarray*}
	\varphi_\e\bigl(\tbeta_1 t,\ldots,\tbeta_d t\bigr)=\sum_{m=1}^{N_\e} q_{y_m} \exp\biggl\{it\sum\limits_{l=1}^d \tr_{m,l} \tbeta_l \biggr\}=\sum_{m=1}^{N_\e} q_{y_m} e^{it y_m}=\poly_{\e}(t),\quad t\in\R,
\end{eqnarray*}
i.e. $\poly_{\e}$ is diagonal for the function $(t_1,\ldots, t_d)\mapsto \varphi_\e(\tbeta_1 t_1,\ldots,\tbeta_d t_d)$. Hence the image of the first function is dense in the image of the second one (see \cite{Levitan}, Theorem 2.4.1., p. 116) and, consequently, in the image of $\varphi_\e$. In particular, this yields
\begin{eqnarray}\label{eq_sup_ge_phie}
	\sup_{t\in\R} |\poly_\e(t)|=\sup_{t\in[-\pi,\pi]^d} |\varphi_\e(t)|.
\end{eqnarray}
Let us show that
\begin{eqnarray}\label{eq_norm_gek_phiek}
	\|\poly_\e^k\|=\|\varphi_\e^k\|\quad\text{for any}\quad k\in\N.
\end{eqnarray}
We fix $k\in\N$ and write
\begin{eqnarray*}
	\poly_\e^k(t)&=&\biggl(\sum_{m=1}^{N_\e} q_{y_m} e^{i ty_m}\biggr)^k=\sum_{m_1=1}^{ N_\e}\ldots\sum_{m_k=1}^{N_\e}\Bigl( q_{y_{m_1}}\cdot\ldots \cdot q_{y_{m_k}} e^{i t(y_{m_1}+\ldots+y_{m_k})}\Bigr),\quad t\in \R,\\
	\varphi_\e^k(t)&=&\biggl(\sum_{m=1}^{N_\e} q_{y_m} e^{i\lng t,\tr_m\rng}\biggr)^k=\sum_{m_1=1}^{ N_\e}\ldots\sum_{m_k=1}^{N_\e}\Bigl( q_{y_{m_1}}\cdot\ldots \cdot q_{y_{m_k}} e^{i\lng t,\tr_{m_1}+\ldots+\tr_{m_k}\rng}\Bigr),\quad t\in \R^d.
\end{eqnarray*}
We next represent
\begin{eqnarray*}
	\poly_\e^k(t)&=&\sum_{z\in \mathcal{Y}_\e^{(k)}}\Biggl(\,\sum_{\substack{m_1=1,\ldots, N_\e,\\ \ldots\\ m_k=1,\ldots, N_\e:\\y_{m_1}+\ldots+y_{m_k}=z}} q_{y_{m_1}}\cdot\ldots \cdot q_{y_{m_k}}\Biggr) e^{itz} ,\quad t\in \R,\\
	\varphi_\e^k(t)&=&\sum_{s\in \mathcal{R}_\e^{(k)}}\Biggl(\,\sum_{\substack{m_1=1,\ldots, N_\e,\\ \ldots\\ m_k=1,\ldots, N_\e:\\ \tr_{m_1}+\ldots+\tr_{m_k}=v}} q_{y_{m_1}}\cdot\ldots \cdot q_{y_{m_k}}\Biggr) e^{i\lng t,s\rng} ,\quad  t\in \R^d,
\end{eqnarray*}
where
\begin{eqnarray*}
	\mathcal{Y}_\e^{(k)}&:=&\Bigl\{y_{m_1}+\ldots+y_{m_k}: m_1,\ldots, m_k\in \{1,\ldots, N_\e\} \Bigr\},\\
	\mathcal{R}^{(k)}_\e&:=&\Bigl\{\tr_{m_1}+\ldots+\tr_{m_k}: m_1,\ldots, m_k\in \{1,\ldots, N_\e\} \Bigr\}.
\end{eqnarray*}
Thus we have
\begin{eqnarray}
	\|\poly_\e^k\|&=&\sum_{z	\in \mathcal{Y}_\e^{(k)}}\Biggl|\sum_{\substack{m_1=1,\ldots, N_\e,\\ \ldots\\ m_k=1,\ldots, N_\e:\\y_{m_1}+\ldots+y_{m_k}=z}} q_{y_{m_1}}\cdot\ldots \cdot q_{y_{m_k}}\Biggr|,\label{eq_norm_gek}\\
	\|\varphi_\e^k\|&=&\sum_{s\in \mathcal{R}_\e^{(k)}}\Biggl|\sum_{\substack{m_1=1,\ldots, N_\e,\\ \ldots\\ m_k=1,\ldots, N_\e:\\ \tr_{m_1}+\ldots+\tr_{m_k}=s}} q_{y_{m_1}}\cdot\ldots \cdot q_{y_{m_k}}\Biggr|.\label{eq_norm_phiek}
\end{eqnarray}
There is a natural map between two finite sets $\mathcal{Y}_\e^{(k)}$ and $\mathcal{R}^{(k)}_\e$: every number $z=\sum_{j=1}^k y_{m_j}$ from the first set is paired with a vector  $s=\sum_{j=1}^k \tr_{m_j}$ from the second set. Here $s$ is the vector of the coefficients of the decomposition for $z$ with the basis $\tbeta_1,\ldots, \tbeta_d$. This map is actually bijection. Indeed, it is injective, because distinct numbers from $\mathcal{Y}_\e^{(k)}$ have the distinct decompositions, i.e. they are paired with the distinct vectors from the set $\mathcal{R}^{(k)}_\e$. The map is surjective, because any vector $v\in \mathcal{R}^{(k)}_\e$ is a sum $\tr_{m'_1}+\ldots+\tr_{m'_k}$ with some indices $m'_1,\ldots, m'_d\in \{1,\ldots,N_\e\}$, and this sum is corresponded to the number $z=y_{m'_1}+\ldots+y_{m'_k}$ from the set $\mathcal{Y}_\e^{(k)}$ by construction. Next, it is clear from the basis decompositions that for any fixed pair of corresponded elements  $z\in \mathcal{Y}_\e^{(k)}$ and $s\in \mathcal{R}^{(k)}_\e$ the equalities $y_{m_1}+\ldots+y_{m_k}=z$ and $\tr_{m_1}+\ldots+\tr_{m_k}=s$ are equivalent for varying index vectors $(m_1,\ldots, m_k)$. Since the inner sums in \eqref{eq_norm_gek} and \eqref{eq_norm_phiek} add the same weight $q_{y_{m_1}}\cdot\ldots \cdot q_{y_{m_k}}$ with every index vector $(m_1,\ldots, m_k)$, we conclude that these inner sums are equal. Thus we come to the equality $\|\poly_\e^k\|=\|\varphi_\e^k\|$. 

We now apply Lemma~\ref{lm_summodulkoefZ} to estimate $\|\varphi_\e^k\|$ for any $k\in\N$ (we set $\varphi:=\varphi_\e$, $N:=N_\e$, $q_m:=q_{y_m}$ and $c_m:=\tr_m$ for every $m=1,\ldots,N_\e$). Namely, using inequality~\eqref{lm_summodulkoefZ_ineq2}, we get 
\begin{eqnarray*}
	\|\varphi_\e^k\|\leqslant A_{\varphi_\e} k^d \sup\limits_{t\in[-\pi,\pi]^d}|\varphi_\e(t)|^k,\quad k\in\N,
\end{eqnarray*}
where $A_{\varphi_\e}$ is a constant defined by \eqref{lm_summodulkoefZ_defAphi}, which doesn't depend on $k$. Applying \eqref{eq_sup_ge_phie} and \eqref{eq_norm_gek_phiek}, we come to the inequality
\begin{eqnarray*}
	\|\poly_\e^k\|\leqslant A_{\varphi_\e} k^d \sup\limits_{t\in\R}|\poly_\e(t)|^k,\quad k\in\N.
\end{eqnarray*}
Observe that
\begin{eqnarray*}
	\sup\limits_{t\in\R}|\poly_\e(t)|\leqslant \sup\limits_{t\in\R}|g_d(t)|+\sup\limits_{t\in\R}|g_d(t)-\poly_\e(t)|.
\end{eqnarray*}
On account of \eqref{ineq_norm_gd-poly}, we have
\begin{eqnarray*}
	\dfrac{\sup_{t\in\R}|g_d(t)-\poly_\e(t)|}{\sup_{t\in\R}|g_d(t)|}\leqslant
	\sup_{t\in\R} \biggl|\dfrac{g_d(t)-\poly_\e(t)}{g_d(t)}\biggr|\leqslant\biggl\|\dfrac{g_d-\poly_\e}{g_d}\biggr\|<\e,
\end{eqnarray*}
i.e. $\sup_{t\in\R}|g_d(t)-\poly_\e(t)|<\e \sup_{t\in\R}|g_d(t)|$. Therefore
\begin{eqnarray*}
	\sup\limits_{t\in\R}|\poly_\e(t)|\leqslant \sup\limits_{t\in\R}|g_d(t)|+\e \sup\limits_{t\in\R}|g_d(t)|=(1+\e) \sup\limits_{t\in\R}|g_d(t)|,
\end{eqnarray*}
where
\begin{eqnarray*}
	\sup\limits_{t\in\R}|g_d(t)|=\sup\limits_{t\in\R}\biggl|\dfrac{1}{f_d(t)}\biggr|=\dfrac{1}{\inf_{t\in\R}|f_d(t)|}=\dfrac{1}{\mu_d}.
\end{eqnarray*}
So we obtain the following estimate
\begin{eqnarray*}
	\|\poly_\e^k\|\leqslant A_{\varphi_\e} k^d\cdot \dfrac{(1+\e)^k}{\mu_d^k},\quad k\in\N.
\end{eqnarray*}
We note that here $k^d\geqslant 1$ and $1/\mu_d^k\geqslant 1$ for any $k\in\N$.  Therefore the estimates of $\|\poly_\e^k\|$  in the cases $Y_\e=\{0\}$ and  $Y_\e\ne\{0\}$ can be unified as follows
\begin{eqnarray*}
	\|\poly_\e^k\|\leqslant C_{\e} k^d\cdot \dfrac{(1+\e)^k}{\mu_d^k}\quad\text{for any} \quad k\in\N,
\end{eqnarray*}
with some constant $C_\e>0$. We use it in \eqref{ineq_gk}:
\begin{eqnarray*}
	\|g_d^k\|\leqslant\dfrac{ \|\poly_\e^k\|}{(1-\e)^k}\leqslant C_{\e} k^d\cdot \dfrac{1}{\mu_d^k}\cdot\biggl(\dfrac{1+\e}{1-\e}\biggr)^k,\quad k\in\N.
\end{eqnarray*}

Let us return to \eqref{ineq_normWn}. Due to the last estimate, for any $\e\in(0,1)$ and $n\in\N$ we have
\begin{eqnarray}\label{ineq_WnVar}
	\|W_n\|\leqslant\sum_{k=1}^{n} \Biggl[ \dfrac{1}{k}\biggl(\dfrac{c_s }{c_d}\biggr)^k  C_{\e} k^d\cdot \dfrac{1}{\mu_d^k}\cdot\biggl(\dfrac{1+\e}{1-\e}\biggr)^k\Biggr]= C_\e \sum_{k=1}^{n} k^{d-1} \biggl(\dfrac{c_s }{c_d\mu_d}\cdot\dfrac{1+\e}{1-\e}\biggr)^k.
\end{eqnarray} 
Since $c_s<c_d \mu_d$ by assumption, the fixed number $\e\in(0,1)$ can be specified as follows
\begin{eqnarray*}
	0<\dfrac{c_s }{c_d\mu_d}\cdot\dfrac{1+\e}{1-\e}<1.
\end{eqnarray*}
Then we obtain
\begin{eqnarray*}
	\sup_{n\in\N}\|W_n\| \leqslant C_{\e}  \sum_{k=1}^{\infty} k^{d-1} \biggl(\dfrac{c_s }{c_d\mu_d}\cdot\dfrac{1+\e}{1-\e}\biggr)^k<\infty,
\end{eqnarray*}
i.e. we showed that the sufficient condition  $2)$ for \eqref{eq_sumintW} is satisfied (the condition $1)$ was proved above). So, according to \eqref{eq_lnfsfd} and \eqref{eq_sumintW}, we have
\begin{eqnarray}\label{eq_inteitxW}
	\ln\biggl(1+ \dfrac{c_s f_s(t)}{c_d f_d(t)}\biggr)=\sum_{k=1}^{\infty} \dfrac{(-1)^{k-1} }{k}\biggl(\dfrac{c_s }{c_d}\biggr)^k f_s(t)^k g_d(t)^k= \int_{\R} e^{itx}\dd W(x),\quad t\in\R,
\end{eqnarray}
with some function $W\in\BV$. Hence  \eqref{eq_Lnfsd} takes the form
\begin{eqnarray*}
	 \Ln f_{s,d}(t)= \int_{\R} e^{itx}\dd W(x)-\ln \biggl(1+ \dfrac{c_s }{c_d}\biggr),\quad t\in\R.
\end{eqnarray*}
The equality $\Ln f_{s,d}(0)=0$, which was mentioned above, implies
\begin{eqnarray*}
\ln \biggl(1+\dfrac{c_s }{c_d}\biggr)= \int_{\R}\dd W(x).
\end{eqnarray*}
Then
\begin{eqnarray*}
	\Ln f_{s,d}(t)= \int_{\R} (e^{itx}-1)\,\dd W(x),\quad t\in\R,
\end{eqnarray*}
i.e. we come to the representation
\begin{eqnarray}\label{repr_fsd}
	f_{s,d}(t)= \exp\biggl\{\int_{\R} (e^{itx}-1)\,\dd W(x)\biggr\}= \exp\biggl\{\int_{\R\setminus\{0\}} (e^{itx}-1)\,\dd W(x)\biggr\},\quad t\in\R.
\end{eqnarray}

Let us investigate the function $W$. If $c_s=0$, then $W$ is identically $0$, which is seen from  \eqref{eq_inteitxW}. We next focus on the case $c_s\ne 0$.  We will prove that $W$ is continuous on $\R$. For it we first observe that the functions $W_n$, $n\in\N$, are continuous on $\R$. Indeed, according to \eqref{def_Wn}, $W_n$ is a finite linear combination of the functions $F_s^{*k}*I_d^{*k}$, $k\in\N$, which are continuous on $\R$ as being  convolutions with continuous $F_s$.  Next, we observe that, similarly to  \eqref{ineq_WnVar},  the following estimate holds
\begin{eqnarray*}
	\|W_{n_2}-W_{n_1}\|\leqslant C_\e\sum_{k=n_1+1}^{n_2} k^{d-1} \biggl(\dfrac{c_s }{c_d\mu_d}\cdot\dfrac{1+\e}{1-\e}\biggr)^k,
\end{eqnarray*}
with the same $\e$, $C_\e$ and  for any positive integers $n_1$ and $n_2$ such that $n_1\leqslant n_2$. The sum from the right-hand side can be  made  arbitrarily small for all sufficiently large $n_1$ and $n_2$, because the terms of this sum (for $k\in\N$) consist the convergent series. This means that $(W_n)_{n\in\N}$ is a fundamental sequence in the space $\BV$. Since $\BV$ is a complete norm space (see the comments at the end of the introduction), there exists $W_*\in\BV$ such that $\|W_n-W_*\|\to 0$ as $n\to\infty$. Hence
\begin{eqnarray*}
	\lim_{n\to\infty}\int_{\R} e^{itx}\dd W_n(x)=\int_{\R} e^{itx}\dd W_*(x),\quad t\in\R.
\end{eqnarray*}
On the other hand, we know that
\begin{eqnarray*}
	\lim_{n\to\infty}\int_{\R} e^{itx}\dd W_n(x)&=&\lim_{n\to\infty}\sum_{k=1}^{n} \dfrac{(-1)^{k-1} }{k}\biggl(\dfrac{c_s }{c_d}\biggr)^k f_s(t)^k g_d(t)^k\\
	&=&\sum_{k=1}^{\infty} \dfrac{(-1)^{k-1} }{k}\biggl(\dfrac{c_s }{c_d}\biggr)^k f_s(t)^k g_d(t)^k=\int_{\R} e^{itx}\dd W(x),\quad t\in\R.
\end{eqnarray*}
Then $\int_{\R} e^{itx}\dd W_*(x)=\int_{\R} e^{itx}\dd W(x)$ for any $t\in\N$, and we conclude that $W_*=W$. So we proved that 
\begin{eqnarray}\label{conv_var_WnW}
	\|W_n-W\|\to 0\quad\text{as}\quad n\to\infty.
\end{eqnarray}
In general, from the well known relation $\sup_{x\in\R}|U(x)|\leqslant \|U\|$ for any $U\in\BV$, we have the uniform convergence:
\begin{eqnarray*}
	\sup_{x\in\R}|W_n(x)-W(x)|\to 0\quad\text{as}\quad n\to\infty.
\end{eqnarray*}
Due to the proved continuity of $W_n$, $n\in\N$, we conclude that the function $W$ is continuous on $\R$ too.

We now show that $W$ can not be purely absolutely continuous function (in the case $c_s\ne 0$). To obtain a contradiction, we suppose that it is not true, i.e. $W$ is absolutely continuous. Let
\begin{eqnarray*}
	w(t):=\int_{\R} e^{itx}\dd W(x),\quad t\in\R.
\end{eqnarray*}
Then, according to \eqref{eq_inteitxW}, we have
\begin{eqnarray*}
	1+ \dfrac{c_s f_s(t)}{c_d f_d(t)}=\exp\biggl\{\int_{\R} e^{itx}\dd W(x)\biggr\}=\exp\{w(t)\}=1+\sum_{k=1}^\infty \dfrac{w(t)^k}{k!},\quad t\in\R.
\end{eqnarray*} 
Hence, on the one hand,
\begin{eqnarray*}
	\dfrac{c_s f_s(t)}{c_d f_d(t)}=\lim_{n\to\infty} \sum_{k=1}^n \dfrac{w(t)^k}{k!}=\lim_{n\to\infty} \sum_{k=1}^n \dfrac{1}{k!} \int_{\R} e^{itx}\dd W^{*k}(x)=\lim_{n\to\infty}  \int_{\R} e^{itx}\dd\biggl(\,\sum_{k=1}^n \dfrac{1}{k!}\,W^{*k}(x)\biggr),\quad t\in\R.
\end{eqnarray*} 
On the other hand, on account of \eqref{eq_fsgd},
\begin{eqnarray*}
	\dfrac{c_s f_s(t)}{c_d f_d(t)}= \dfrac{c_s }{c_d}\cdot f_s(t)\cdot g_d(t)=\dfrac{c_s }{c_d}\cdot \int_{\R} e^{itx}\dd (F_s*I_d)(x)= \int_{\R} e^{itx}\dd \biggl(\dfrac{c_s }{c_d}\cdot (F_s*I_d)(x)\biggr),\quad t\in\R.
\end{eqnarray*} 
Since the variance
\begin{eqnarray*}
	\biggl\|\sum_{k=n_1+1}^{n_2} \dfrac{1}{k!}\,W^{*k}\biggr\|\leqslant \sum_{k=n_1+1}^{n_2} \dfrac{\|W^{*k}\|}{k!}\leqslant \sum_{k=n_1+1}^{n_2} \dfrac{\|W\|^k}{k!}
\end{eqnarray*}
can be  made  arbitrarily small for all sufficiently large $n_1$ and $n_2$, the sequence of sums $\sum_{k=1}^n \tfrac{1}{k!}\,W^{*k}$, $n\in\N$, is fundamental in the space $\BV$. Due to the completeness of $\BV$,  these sums converge in variation to  some function from $\BV$, namely, to $\tfrac{c_s }{c_d}\cdot (F_s*I_d)$ by the uniqueness of the Fourier--Stieltjes transform, i.e.
\begin{eqnarray*}
	\biggl\|\sum_{k=1}^{n} \dfrac{1}{k!}\,W^{*k}- \dfrac{c_s }{c_d}\cdot (F_s*I_d)\biggr\|\to 0\quad\text{as} \quad n\to\infty.
\end{eqnarray*}
But this is impossible, because $\sum_{k=1}^{n} \tfrac{1}{k!}\,W^{*k}$, $n\in\N$, are absolutely continuous as linear combinations of convolution powers of absolutely continuous $W$ by the assumption, but $\tfrac{c_s }{c_d}\cdot (F_s*I_d)$ is continuous singular as the convolution of the continuous singular function $F_s$ and the discrete function $I_d$ (see  the comments below Remark~\ref{rem_notdomin}). Thus $W$ is not (purely) absolutely continuous function from $\BV$, i.e. it always has some continuous singular part.

Let us return to \eqref{conv_var_WnW} and prove the sufficient condition (from the statement of the theorem) for $W$ to be (purely) continuous singular. Suppose that all the functions $F_s^{*k}$, $k\in\N$,  are continuous singular. Hence all $W_n$, $n\in\N$ are too. Let  $W_a$ and $W_s$ be the absolutely continuous part and the  continuous singular part of the Lebesgue decomposition for $W$, respectively, i.e. $W=W_a+W_s$. Then
\begin{eqnarray*}
	\|W_n-W\|= \|W_n-W_a-W_s\|=\|W_n-W_s\|+\|W_a\|\geqslant \|W_a\|\geqslant 0.
\end{eqnarray*}
Due to \eqref{conv_var_WnW}, we conclude that $\|W_a\|=0$, i.e. $W_a(x)=0$ for any $x\in\R$. Thus $W=W_s$.  

We now combine representations \eqref{repr_fd}, \eqref{repr_fads} and \eqref{repr_fsd} in formula \eqref{repr_fd_fads_fsd}:
\begin{eqnarray*}
	f(t)&=&\exp\Biggl\{ i t\gamma_0  + \sum_{u \in \langle \X \rangle \setminus\{0\}} \lambda_u \bigr(e^{i tu} - 1\bigr)\\
	 &&{}\qquad\qquad+\int_{\R\setminus\{0\}} (e^{itx}-1)  \biggl( v_a(x)+\sgn(x)\,  \dfrac{\mathfrak{m}_a\cdot e^{-|x|}}{|x|}\biggr) \dd x+\int_{\R\setminus\{0\}} (e^{itx}-1)\,\dd W(x)\Biggr\}\\
	 &=&\exp\biggl\{ i t\gamma_0  + \int_{\R\setminus\{0\}} (e^{itx}-1)\,\dd L(x)\biggr\}
	 ,\quad t\in\R,
\end{eqnarray*}
with $L(x):=L_d(x)+L_a(x)+L_s(x)$, $x\in\R\setminus\{0\}$, where $L_d$, $L_a$, and $L_s$ are defined by formulas \eqref{def_Ld}, \eqref{def_La}, and \eqref{def_Ls}, respectively. It is not difficult to check that  $L_d$, $L_a$, $L_s$, and, consequently, $L$ satisfy all conditions for spectral function of the L\'evy type representation (see introduction or Theorem~\ref{th_Berger2}). However, we recall that the function $v_a$ is potentially complex-valued and hence  $L$ is too.  Next, it is seen that
\begin{eqnarray*}
	\int_{S_1}|x|\,\dd |L_*|(x)<\infty
\end{eqnarray*}
with $S_1:=[-1,1]\setminus\{0\}$ and for any $L_*\in\{L_d,L_a,L_s\}$. Hence it is also true for $L_*=L$, and we can write
\begin{eqnarray*}
	f(t)=\exp\biggl\{ i t\gamma_1  + \int_{\R\setminus\{0\}} \Bigl(e^{itx}-1-it x\,\idd_{[-1,1]}(x)\Bigr)\,\dd L(x)\biggr\}
	,\quad t\in\R,
\end{eqnarray*}
with $\gamma_1:=\gamma_0+ \int_{\R\setminus\{0\}} x\,\idd_{[-1,1]}(x)\,\dd L(x)$, which is potentially complex.  

We now apply Theorem~\ref{th_Berger2}. Accordingly, $L$ is actually real-valued, which means that $v_a$ is a real-valued function. Moreover, we have $F\in\RID$, i.e. $(i)$ is proved.\quad $\Box$\\

\textbf{Proof of Proposition~\ref{pr_PureSing}.} To obtain a contradiction, we assume that $W$ is continuous singular.  
For convenience, we introduce the operator $[\,\cdot\,]_s$, which acts from $\BV$ to $\BV$, and for any  $U\in\BV$ it returns the continuous singular part of $U$ as $[U]_s$ (we set $[U]_s$ to be identically zero if there is no continuous singular part). We start with \eqref{conv_var_WnW} from the proof of Theorem~\ref{th_MainResult}. Under the assumptions of the proposition, it holds. For every $n\in\N$ we write
\begin{eqnarray*}
	\|W_n-W\|= \bigl\|(W_n-[W_n]_s)+([W_n]_s-W)\bigr\|=\bigl\|W_n-[W_n]_s\bigr\|+\bigl\|[W_n]_s-W\bigr\|,
\end{eqnarray*}
where the latter equality is valid, because $W_n-[W_n]_s$ is absolutely continuous (or identically zero) and $[W_n]_s-W$ is continuous singular by the assumption. Therefore $\|W_n-W\|\geqslant \bigl\|[W_n]_s-W\bigr\|\geqslant 0$ and, due to \eqref{conv_var_WnW}, we conclude that
\begin{eqnarray}\label{conv_WnSingW}
	\bigl\|[W_n]_s-W\bigr\|\to 0\quad\text{as}\quad n\to\infty.
\end{eqnarray}
Next, according to \eqref{def_Wn}, we have
\begin{eqnarray*}
	[W_n]_s=\Biggl[\,\sum_{k=1}^{n}  \dfrac{(-1)^{k-1} }{k}\biggl(\dfrac{c_s }{c_d}\biggr)^k  F_s^{*k}*I_d^{*k}\Biggr]_s= \sum_{k=1}^{n}  \dfrac{(-1)^{k-1} }{k}\biggl(\dfrac{c_s }{c_d}\biggr)^k  \bigl[F_s^{*k}*I_d^{*k}\bigr]_s,\quad n\in\N.
\end{eqnarray*} 
Since $I_d$ is discrete (see \eqref{def_Id}),  the functions $I_d^{*k}$ are discrete and hence $\bigl[F_s^{*k}*I_d^{*k}\bigr]_s=\bigl[F_s^{*k}\bigr]_s*I_d^{*k}$ for any $k\in\N$. Thus we come to the equalities
\begin{eqnarray}\label{eq_WnSing}
	[W_n]_s= \sum_{k=1}^{n}  \dfrac{(-1)^{k-1} }{k}\biggl(\dfrac{c_s }{c_d}\biggr)^k  \bigl[F_s^{*k}\bigr]_s*I_d^{*k},\quad n\in\N.
\end{eqnarray} 
Let us consider the functions $\bigl[F_s^{*k}\bigr]_s$ and numbers $\bigl[F_s^{*k}\bigr]_s(\infty)$. By the definition of $n_a$, for any $k<n_a$ we have $\bigl[F_s^{*k}\bigr]_s=F_s^{*k}$ and, in particular, $\bigl[F_s^{*k}\bigr]_s(\infty)=1$. For the case $k=n_a$, by the assumption, we have  $\bigl[F_s^{*n_a}\bigr]_s(x)= (1-\alpha)H_s(x)$, $x\in\R$, and hence $\bigl[F_s^{*n_a}\bigr]_s(\infty)= 1-\alpha$. Next, for any integer $k\geqslant n_a$ we observe that
\begin{eqnarray*}
	\bigl[F_s^{*(k+1)}\bigr]_s(x)=\bigl[F_s^{*k}*F_s\bigr]_s(x)=\bigl[\bigl[F_s^{*k}\bigr]_s*F_s\bigr]_s(x)\leqslant \bigl([F_s^{*k}\bigr]_s*F_s\bigr)(x),\quad x\in\R.
\end{eqnarray*}
By the Lebesgue dominated convergence theorem, we conclude that
\begin{eqnarray*}
	\bigl([F_s^{*k}\bigr]_s*F_s\bigr)(\infty)=\lim\limits_{z\to\infty} \int_{\R} F_s(z-x)\, \dd [F_s^{*k}\bigr]_s(x)=\int_{\R} \dd [F_s^{*k}\bigr]_s(x)=\bigl[F_s^{*k}\bigr]_s(\infty).
\end{eqnarray*}
Therefore  $\bigl[F_s^{*(k+1)}\bigr]_s(\infty)\leqslant\bigl[F_s^{*k}\bigr]_s(\infty)$ for any integer $k\geqslant n_a$. Let us introduce the sequence $A_k:= \bigl[F_s^{*k}\bigr]_s(\infty)$, $k\in\N$. So we have
\begin{eqnarray}\label{eq_Ak}
	A_1=\ldots=A_{n_a-1}=1,\qquad A_{n_a}=1-\alpha, \qquad A_{k}\geqslant A_{k+1}\geqslant 0 \quad\text{for any}\quad k\geqslant n_a.
\end{eqnarray}

Due to \eqref{conv_WnSingW}, we next conclude that
\begin{eqnarray*}
	\int_{\R} e^{itx}\dd W(x)=\lim_{n\to\infty}\int_{\R} e^{itx}\dd [W_n]_s(x),\quad t\in\R.
\end{eqnarray*}
Then, according to \eqref{eq_inteitxW} and \eqref{eq_WnSing},  we have
\begin{eqnarray*}
	\ln\biggl(1+ \dfrac{c_s f_s(t)}{c_d f_d(t)}\biggr)=\lim_{n\to\infty}\sum_{k=1}^{n} \dfrac{(-1)^{k-1} }{k}\biggl(\dfrac{c_s }{c_d}\biggr)^k f_{k,s}(t) g_d(t)^k \quad\text{for any} \quad t\in\R,
\end{eqnarray*}
where
\begin{eqnarray*}
	f_{k,s}(t):= \int_{\R} e^{itx}\,\dd\bigl[F_s^{*k}\bigr]_s(x),\quad t\in\R.
\end{eqnarray*}
In particular, we write
\begin{eqnarray*}
	\ln\biggl(1+ \dfrac{c_s f_s(0)}{c_d f_d(0)}\biggr)=\sum_{k=1}^{\infty} \dfrac{(-1)^{k-1} }{k}\biggl(\dfrac{c_s }{c_d}\biggr)^k f_{k,s}(0) g_d(0)^k.
\end{eqnarray*}
Recall that $f_s(0)=f_d(0)=g_d(0)=1$ and $f_{k,s}(0)= [F_s^{*k}\bigr]_s(\infty)=A_k$ for any $k\in\N$. Then, on the one hand,
\begin{eqnarray*}
	\ln\biggl(1+ \dfrac{c_s}{c_d }\biggr)=\sum_{k=1}^{\infty} \dfrac{(-1)^{k-1} }{k}\biggl(\dfrac{c_s }{c_d}\biggr)^k A_k.
\end{eqnarray*}
On the other hand, since $c_s<c_d\mu_d\leqslant c_d$, we have the expansion
\begin{eqnarray*}
	\ln\biggl(1+ \dfrac{c_s}{c_d }\biggr)=\sum_{k=1}^{\infty} \dfrac{(-1)^{k-1} }{k}\biggl(\dfrac{c_s }{c_d}\biggr)^k.
\end{eqnarray*}
Thus we come to the following equality
\begin{eqnarray*}
	\sum_{k=1}^{\infty} \dfrac{(-1)^{k-1} }{k}\biggl(\dfrac{c_s }{c_d}\biggr)^k A_k=\sum_{k=1}^{\infty} \dfrac{(-1)^{k-1} }{k}\biggl(\dfrac{c_s }{c_d}\biggr)^k.
\end{eqnarray*}
Due to \eqref{eq_Ak}, it is reduced to
\begin{eqnarray}\label{eq_series}
	\sum_{k=n_a}^{\infty} \dfrac{(-1)^{k-n_a} }{k}\biggl(\dfrac{c_s }{c_d}\biggr)^k A_k=\sum_{k=n_a}^{\infty} \dfrac{(-1)^{k-n_a} }{k}\biggl(\dfrac{c_s }{c_d}\biggr)^k.
\end{eqnarray}
The series from the left-hand side is alternating with non-increasing absolute values of the terms. Then the following inequality is valid:
\begin{eqnarray*}
	\sum_{k=n_a}^{\infty} \dfrac{(-1)^{k-n_a} }{k}\biggl(\dfrac{c_s }{c_d}\biggr)^k A_k\leqslant  \dfrac{1 }{n_a}\biggl(\dfrac{c_s }{c_d}\biggr)^{n_a} A_{n_a}=\dfrac{1 }{n_a}\biggl(\dfrac{c_s }{c_d}\biggr)^{n_a}(1-\alpha).
\end{eqnarray*}
Since the series from the right-hand side of \eqref{eq_series} is alternating with strictly decreasing absolute values of the terms,  the following inequality hold:
\begin{eqnarray*}
	\sum_{k=n_a}^{\infty} \dfrac{(-1)^{k-n_a} }{k}\biggl(\dfrac{c_s }{c_d}\biggr)^k>\dfrac{1 }{n_a}\biggl(\dfrac{c_s }{c_d}\biggr)^{n_a}-\dfrac{1 }{n_a+1}\biggl(\dfrac{c_s }{c_d}\biggr)^{n_a+1}=\dfrac{1 }{n_a}\biggl(\dfrac{c_s }{c_d}\biggr)^{n_a}\biggl(1-\dfrac{n_a }{n_a+1}\cdot\dfrac{c_s }{c_d}\biggr).
\end{eqnarray*}
We see that the assumption $\alpha\geqslant \tfrac{n_a }{n_a+1}\cdot\tfrac{c_s }{c_d}$ implies the inequality
\begin{eqnarray*}
	\sum_{k=n_a}^{\infty} \dfrac{(-1)^{k-n_a} }{k}\biggl(\dfrac{c_s }{c_d}\biggr)^k A_k<\sum_{k=n_a}^{\infty} \dfrac{(-1)^{k-n_a} }{k}\biggl(\dfrac{c_s }{c_d}\biggr)^k.
\end{eqnarray*}
Thus we come  to a contradiction with \eqref{eq_series}. So the assumption, which we started with, is false that $W$ is continuous singular.

Due to the condition of the dominated continuous singular part, we have $c_s<c_d$. Since $\tfrac{n_a}{n_a+1}<1$, the inequality $\alpha\geqslant \tfrac{n_a }{n_a+1}\cdot\tfrac{c_s }{c_d}$ always holds in the case $\alpha=1$.\quad $\Box$\\

\textbf{Proof of Proposition \ref{pr_Example}.} Let us consider the function $f_s$ defined by formula \eqref{def_fs}.
We first find the limit of the sequence $f_s(t_n)$, $n\in\N$, with $t_n:=\pi(2n)!$ as $n\to\infty$. We write $f_s(t_n)=M_n\cdot R_n$, $n\in\N$, where
\begin{eqnarray*}
	M_n:=\prod\limits_{k=1}^{2n} \cos(t_n/k!)\quad\text{and}\quad R_n:=\prod\limits_{k=2n+1}^{\infty} \cos(t_n/k!),\quad n\in\N.
\end{eqnarray*}
For any $n\in\N$ 
\begin{eqnarray*}
	M_n&=&\cos \biggl(\dfrac{\pi (2n)!}{(2n)!}\biggr)\cdot \cos \biggl(\dfrac{\pi (2n)!}{(2n-1)!}\biggr)\cdot \cos \biggl(\dfrac{\pi (2n)!}{(2n-2)!}\biggr)\cdot\ldots \cdot \cos \biggl(\dfrac{\pi (2n)!}{1!}\biggr)\\
	&=&\cos(\pi)\cdot \cos(\pi\cdot2n)\cdot \cos\bigl(\pi\cdot 2n(2n-1)\bigr)\cdot\ldots\cdot \cos\bigl(\pi\cdot 2n(2n-1)\cdot\ldots\cdot 1\bigr)\\
	&=&(-1)\cdot 1\cdot 1\cdot\ldots \cdot 1=-1.  
\end{eqnarray*}
Next, for every $n\in\N$ and for every integer $k\geqslant 2n+1$ we have
\begin{eqnarray*}
	0<\dfrac{t_n}{k!}= \dfrac{\pi(2n)!}{k\cdot(k-1)!}\leqslant \dfrac{\pi}{k}\leqslant \dfrac{\pi}{3},
\end{eqnarray*}
and, by the well-known inequality $\cos(x)\geqslant 1-x^2/2$, $x\in\R$, we get
\begin{eqnarray*}
	\cos(t_n/k!)\geqslant 1-\dfrac{(t_n/k!)^2}{2}\geqslant 1-\dfrac{\pi^2}{2k^2}\geqslant 1-\dfrac{\pi^2}{18}>0.
\end{eqnarray*}
Therefore, on the one hand, it is clear that for any $n\in\N$ we have $\cos(t_n/k!) <1$ for any $k\geqslant 2n+1$ and hence $R_n< 1$. On the other hand,
\begin{eqnarray*}
	R_n\geqslant\prod\limits_{k=2n+1}^{\infty} \biggl(1-\dfrac{\pi^2}{2k^2}\biggr)=\exp\biggl\{\sum_{k=2n+1}^\infty \ln\biggl(1-\dfrac{\pi^2}{2k^2}\biggr) \biggr\},\quad n\in\N,
\end{eqnarray*}
where the sum in the exponent tends to $0$ as $n\to\infty$. So we conclude that $R_n\to 1$ as $n\to\infty$. Thus
\begin{eqnarray*}
	f_s(t_n)=M_n\cdot R_n=-R_n\to-1 \quad\text{as}\quad n\to\infty,
\end{eqnarray*}
and we know that $f_s(t_n)>-1$ for any $n\in\N$.

We next observe that
\begin{eqnarray*}
	f_s(t_n\pm \pi)=\prod\limits_{k=1}^{\infty} \cos\biggl(\dfrac{t_n\pm \pi}{k!}\biggr)=0.
\end{eqnarray*}
Indeed, 
\begin{eqnarray*}
	\cos\biggl(\dfrac{t_n\pm \pi}{k!}\biggr)\bigg|_{k=2}=\cos\biggl(\dfrac{\pi (2n)!\pm \pi}{2!}\biggr)=\cos\biggl(\pi \cdot 3\cdot\ldots\cdot (2n)\pm \dfrac{\pi}{2}\biggr)=0.
\end{eqnarray*}

We now return to the characteristic function $f_*(t)=\tfrac{1}{2}+\tfrac{1}{2}\,f_s(t)$, $t\in\R$, and we consider the quantities
\begin{eqnarray*}
	\dfrac{f_*(t_n-\pi)f_*(t_n+\pi)}{f_*(t_n)^2}=\dfrac{ \bigl(\tfrac{1}{2}+\tfrac{1}{2}\,f_s(t_n-\pi)\bigr)\bigl(\tfrac{1}{2}+\tfrac{1}{2}\,f_s(t_n+\pi)\bigr) }{\bigl(\tfrac{1}{2}+\tfrac{1}{2}\,f_s(t_n)\bigr)^2},\quad n\in\N.
\end{eqnarray*}
By the above, 
\begin{eqnarray*}
	\dfrac{f_*(t_n-\pi)f_*(t_n+\pi)}{f_*(t_n)^2}=\dfrac{1}{\bigl(1+f_s(t_n)\bigr)^2}\to\infty\quad\text{as}\quad n\to\infty.
\end{eqnarray*}
Therefore, by  Theorem~\ref{th_psifunction}, $f_*$ can not be characteristic function for a distribution function from the class $\RID$. Thus $F_*\notin \RID$.\quad $\Box$\\

\textbf{Proof of Proposition \ref{pr_Decom}.} We have decompositions \eqref{repr_F_Lebesgue} and \eqref{repr_f_Lebesgue} for the distribution function $F$ and its characteristic function $f$, respectively. Since $F\in\RID$, we know that $f(t)\ne 0$ for any $t\in\R$ and $\mu_d=\inf_{t\in\R}|f_d(t)|>0$ according to Theorem \ref{th_MainResult_Levy}. Due to the assumption of the dominated singular part and $c_d>0$, we have $c_s<c_d\mu_d$.

Let $f_1$ and $f_2$ denote correspondingly the characteristic functions of $F_1$ and $F_2$. Then the assumed decomposition  $F=F_1*F_2$ means that $f(t)=f_1(t)f_2(t)$, $t\in\R$. Hence we immediately note that $f_1(t)\ne 0$  and $f_2(t)\ne 0$ for any $t\in\R$. 

Let us write the Lebesgue decompositions for $F_1$ and $F_2$:
\begin{eqnarray*}
	F_j(x)=c_{j,d}F_{j,d}(x)+c_{j,a}F_{j,a}(x)+c_{j,s}F_{j,s}(x),\quad x\in\R,\quad j\in\{1,2\}.
\end{eqnarray*}
We also write the corresponding decompositions for $f_1$ and $f_2$:
\begin{eqnarray*}
	f_j(t)=c_{j,d}f_{j,d}(t)+c_{j,a}f_{j,a}(t)+c_{j,s}f_{j,s}(t),\quad t\in\R,\quad j\in\{1,2\}.
\end{eqnarray*}
Here $c_{j,d}$, $c_{j,a}$, $c_{j,s}$ are non-negative and $c_{j,d}+c_{j,a}+c_{j,s}=1$ for  $j=1$ and $j=2$.
For clarity, we write the equality $F=F_1*F_2$ in the expanded form:
\begin{eqnarray}\label{eq_decomFF1F2}
	c_{d}F_{d}+c_{a}F_{a}+c_{s}F_{s}=\bigl(c_{1,d}F_{1,d}+c_{1,a}F_{1,a}+c_{1,s}F_{1,s}\bigr)*\bigl(c_{2,d}F_{2,d}+c_{2,a}F_{2,a}+c_{2,s}F_{2,s}\bigr).
\end{eqnarray}
Since $F$ has non-zero discrete part ($c_d>0$), the functions $F_1$ and $F_2$  have non-zero discrete parts too, i.e. $c_{1,d}>0$ and $c_{2,d}>0$. Since a convolution of any two distribution functions is discrete if and only if  these functions are discrete, we conclude that $c_dF_d(x)=c_{1,d}c_{2,d}(F_{1,d}*F_{2,d})(x)$, $x\in\R$, i.e. $c_d =c_{1,d}c_{2,d}$ and $F_d=F_{1,d}*F_{2,d}$. Thus we have $f_d(t)=f_{1,d}(t)f_{2,d}(t)$, $t\in\R$. Since $|f_{1,d}(t)|\leqslant 1$ and $|f_{2,d}(t)|\leqslant 1$ for any $t\in\R$, we conclude that
\begin{eqnarray*}
	\mu_{1,d}:=\inf_{t \in \R}|f_{1,d}(t)|\geqslant \inf_{t \in \R}|f_{d}(t)|=\mu_d>0,
\end{eqnarray*}
and, analogously,
\begin{eqnarray*}
	\mu_{2,d}:=\inf_{t \in \R}|f_{2,d}(t)|\geqslant \inf_{t \in \R}|f_{d}(t)|=\mu_d>0.
\end{eqnarray*}

We next observe that $F_{1,s}*F_{2,d}$ and  $F_{2,s}*F_{1,d}$ are continuous singular. Therefore the corresponding summands from \eqref{eq_decomFF1F2}  are included in the continuous part of $F$, i.e. $c_sF_s(x)\geqslant c_{1,s}c_{2,d}(F_{1,s}*F_{2,d})(x) +c_{2,s}c_{1,d}(F_{2,s}*F_{1,d})(x)$, $x\in\R$. Consequently,
\begin{eqnarray*}
	c_s=c_s F_s(\infty)\geqslant  c_{1,s}c_{2,d}(F_{1,s}*F_{2,d})(\infty) +c_{2,s}c_{1,d}(F_{2,s}*F_{1,d})(\infty)  =c_{1,s}c_{2,d}+c_{2,s}c_{1,d }.
\end{eqnarray*}
Then we get
\begin{eqnarray*}
	c_{1,s}c_{2,d}\leqslant c_s<c_d\mu_d=c_{1,d} c_{2,d}\mu_d\leqslant c_{1,d} c_{2,d}\mu_{1,d},
\end{eqnarray*}
i.e. $c_{1,s}<c_{1,d}\mu_{1,d}$. Analogously,
\begin{eqnarray*}
	c_{2,s}c_{1,d}\leqslant c_s<c_d\mu_d=c_{1,d} c_{2,d}\mu_d\leqslant c_{1,d} c_{2,d}\mu_{2,d},
\end{eqnarray*}
i.e. $c_{2,s}<c_{2,d}\mu_{2,d}$. Thus $F_1$ and $F_2$ have the dominated continuous singular parts.

So we have showed that $F_1$ and $F_2$ satisfy the assumptions of Theorem~\ref{th_MainResult} and the condition $(iii)$ from it. Thus, according to this theorem,  $F_1$ and $F_2$ belong to the class $\RID$.  \quad $\Box$


\begin{thebibliography}{99}
	\bibitem{AlexeevKhartov}  I.\,A. Alexeev, A.\, A. Khartov, \textit{Spectral representations of characteristic functions of discrete probability laws}, Bernoulli \textbf{29} (2023), 2, 1392--1409.	
	
 	\bibitem{AlexeevKhartov_Mult} I.\,A. Alexeev, A.\,A. Khartov, \textit{A criterion and a Crame\'er-Wold device for quasi-infinite divisibility for discrete multivariate probability laws}, Electron. J. Probab., \textbf{28} (2023), 1--17.	
	
	\bibitem{AlexeevKhartov2} I. A. Alexeev, A. A. Khartov, \textit{On convergence and compactness in variation with a shift of discrete probability laws}, Vestnik St.Petersburg University: Mathematics, \textbf{54} (2021), issue 3, p. 221--226.	
		
	\bibitem{Appl} D. Applebaum, \textit{L\'evy Processes and Stochastic Calculus}, Camb. Univ. Press, Cambridge, 2009.	
		
	\bibitem{BaxterShapiro} G. Baxter, J.\,M. Shapiro, \textit{On bounded infinitely divisible random variables}, Sankhya: The Indian Journal of Statistics, \textbf{22}  (1960), No. 3/4, 253--260.
	
	\bibitem{Berger}  D. Berger, \textit{On quasi-infinitely divisible distributions with a point mass}, Math. Nachr., \textbf{292} (2019), 1674--1684.
	
	\bibitem{BergerKutlu} D. Berger, M. Kutlu, \textit{Quasi-infinite divisibility of a class of distributions with discrete part}, Proc. Amer. Math. Soc. \textbf{151} (2023), 5, 2211--2224.
	
	\bibitem{BergKutLind} D. Berger, M. Kutlu, A. Linder, \textit{On multivariate quasi-infinitely divisible distributions}, A Lifetime of Excersions Through Random Walks and L\'evy Processes. A Volume in Honour of Ron Doney's 80th Birthday. L. Chaumont, A.E. Kyprianou (eds.), Progress in Probability 78 (2021), Birkh\"auser,  87--120.	
		
	\bibitem{Bochner} S. Bochner, A theorem on Fourier-Stieltjes integrals, Bull. Amer. Math. Soc., \textbf{40} (1934), 4, 271--276.	
		
		
	\bibitem{ChhDemniMou}   H. Chhaiba, N. Demni, Z.  Mouayn, \textit{Analysis of generalized negative binomial distributions attached to hyperbolic Landau levels}, J. Math. Phys., \textbf{57} (2016), 7, 072103.	
		
	\bibitem{Gelfand} I.\,M. Gelfand, D.\,A. Raikov, G.\,E. Shilov, \textit{Commutative Normed Rings}, New York: Chelsea Publishing Company, 1964.
		
	\bibitem{GnedKolm} B.\,V. Gnedenko,  A.\,N. Kolmogorov,  \textit{Limit Distributions for Sums of Independent Random Variables}, Cambridge: Addison-Wesley, 1954.	
	
	\bibitem{Hardy} M. Hardy, \textit{Combinatorics of Partial Derivatives}, Electron. J. Comb., 13 (2006), R1.
	
	\bibitem{HewittZuckermann} E. Hewitt, H.\,S. Zuckerman, \textit{Singular measures with absolutely continuous convolution squares}, Proc. Camb. Phil. Soc., \textbf{62} (1966), , 399--420.
		
	\bibitem{IvashevMusatov} O.\,S. Iva$\check{\mathrm{s}}$ev-Musatov,  \textit{On coefficients of trigonometric null-series}. Izv. Akad. Nauk SSSR, Ser. Mat., \textbf{21} (1957), 559--578.	
		
		
	\bibitem{JessenWintner} B. Jessen, A. Wintner, \textit{Distribution functions and the Riemann zeta function}, Trans. Amer. Math. Soc., \textbf{38} (1935), 48--88.
		
	 \bibitem{KhartovComp} A.\,A. Khartov, \textit{Compactness criteria for quasi-infinitely divisible distributions on the integers}, Stat. Probab. Lett., \textbf{153} (2019), 1--6.	
		
	\bibitem{Khartov} A.\,A. Khartov, \textit{A criterion of quasi-infinite divisibility for discrete laws}, Stat. Probab. Lett., \textbf{185} (2022), 109436.	
		
	\bibitem{KhartovWeak}	A.\,A. Khartov, \textit{On weak convergence of quasi-infinitely divisible laws}, Pacific Journal of Mathematics, \textbf{322} (2023), 2, 341--367.	
		

	\bibitem{KhartovGeneral} A.\,A. Khartov, \textit{Some criteria of rational-infinite divisibility for probability laws}, arXiv:2305.14524v2.
	
	\bibitem{KhartovDecomQ} A.\,A. Khartov, \textit{On decomposition problem for distribution functions of class $\RID$}, arXiv:2412.18915.
	
	
	\bibitem{KhartovAlexeevThree} A.\,A. Khartov, I.\,A. Alexeev,  \textit{Quasi-infinite divisibility and three-point probability laws}, J. Math. Sci. \textbf{268} (2022), 731--738.
	
	\bibitem{Krein} M.\,G. Krein, \textit{Integral equations on the half-line with a kernel depending on the difference of the arguments}, Uspekhi Mat. Nauk, \textbf{13} (1958), 5, 3--120 (in Russian). 
	

	
	
	
	\bibitem{Levitan} B. M. Levitan, \textit{Almost periodic functions}, GITTL, Moskow, 1953 (In Russian).
	
	\bibitem{Lifshits}  M.\,A. Lifshits, \textit{Random Processes by Example}, World Scientific, Singapore, 2014.
	
	\bibitem{LindPanSato} A. Lindner, L. Pan, K. Sato, \textit{On quasi-infinitely divisible distributions}, Trans. Amer. Math. Soc., \textbf{370} (2018), 8483--8520.
	
	\bibitem{LinOstr} Yu.\,V. Linnik, I.\,V. Ostrovskii, \textit{Decomposition of Random Variables and Vectors}, Transl. Math. Monog. \textbf{48}, AMS, Providence, Rhode Island, 1977.
	
	\bibitem{Lukacs} E. Lukacs, \textit{Characteristic Functions}, Griffin, London, 1970.
	
	\bibitem{Nakamura} T. Nakamura,  \textit{A complete Riemann zeta distribution and the Riemann hypothesis}, Bernoulli \textbf{21} (2015), 1, 604--617.
	
	\bibitem{Passegg} R. Passeggeri, \textit{On Quasi-Infinitely Divisible Random Measures}. Bayesian Anal., \textbf{18} (2023), 1, 253--286.
	
	
	\bibitem{Sato1999}  K.-I. Sato, \textit{L\'evy Processes and Infinitely Divisible Distributions}, Camb. Univ. Press, Cambridge, 1999.
	
	\bibitem{Schoutens} W. Schoutens, \textit{L\'evy Processes in Finance: Pricing Financial Derivatives}, Wiley, Chichester, 2003.
	
	\bibitem{Shreider} Yu.\,A. Shreider, \textit{The structure of maximal ideals in rings of measures with convolution}, Sbornik Math., \textbf{69} (1950), 2, 297--318 (in Russian).
	
	
	
	\bibitem{Tucker} H.\,G. Tucker, \textit{On a necessary and sufficient condition that an infinitely divisible distribution be absolutely continuous}, Trans. Amer. Math. Soc., \textbf{118} (1965), 316--330.
	
	
	\bibitem{WienerPitt} N. Wiener, H.\,R. Pitt, \textit{On absolutely convergent Fourier-Stieltjes transforms}, Duke Math. J., \textbf{4} (1938), 420--430.
	
	\bibitem{WienerWintner} N. Wiener, A. Wintner, Fourier-Stieltjes transforms and singular infinite convolutions, Amer. J. Math., \textbf{60} (1938), 3, 513--522.
	
	
	\bibitem{Zolot2} V.\,M. Zolotarev,  \textit{Modern Theory of Summation of Random Variables},  Utrecht: VSP, 1997.
	\end{thebibliography}
\end{document}